\documentclass[11pt,reqno]{amsart}

\usepackage{amsmath, amssymb, amsthm}
\usepackage{mathrsfs}
\usepackage{mathtools}
\usepackage{setspace}

\usepackage{hyperref}
\hypersetup{colorlinks=true, citecolor=blue, linkcolor=blue, urlcolor=blue, pdfstartview=FitH}
\usepackage{cleveref}

\usepackage[dvipsnames]{xcolor}
\usepackage{comment}

\usepackage{tikz}
\usepackage{tikz-cd}
\usepackage[all,2cell,cmtip]{xy}
\UseTwocells

\theoremstyle{plain} %
\newtheorem{theorem}{Theorem}[section]
\newtheorem*{NAHC}{Fiberwise Nonabelian Hodge Correspondence}
\newtheorem*{RNAHC}{Universal Nonabelian Hodge Correspondence}
\newtheorem{prop}[theorem]{Proposition}
\newtheorem{lemma}[theorem]{Lemma}
\newtheorem{corollary}[theorem]{Corollary}
\newtheorem{definition}[theorem]{Definition}

\theoremstyle{definition}
\newtheorem{remark}[theorem]{Remark}

\newcommand{\QQ}{\mathbb{Q}}
\newcommand{\CC}{\mathbb{C}}
\newcommand{\PP}{\mathbb{P}}
\newcommand{\Hom}{\mathrm{Hom}}
\newcommand{\M}{\mathcal{M}}
\newcommand{\C}{\mathcal{C}}
\newcommand{\mc}{\mathcal}

\newcommand{\ZZ}{\mathbb{Z}}
\DeclareMathOperator{\Id}{Id}
\DeclareMathOperator{\Diff}{Diff}
\DeclareMathOperator{\Sp}{Sp}
\DeclareMathOperator{\GL}{GL}
\DeclareMathOperator{\SL}{SL}
\DeclareMathOperator{\PGL}{PGL}
\newcommand{\BGL}{B\mathrm{GL}}
\newcommand{\BSL}{B\mathrm{SL}}
\newcommand{\BPGL}{B\mathrm{PGL}}
\newcommand{\BSO}{B\mathrm{SO}}
\DeclareMathOperator{\Mod}{Mod}
\DeclareMathOperator{\Tr}{tr}
\DeclareMathOperator{\Sym}{Sym}
\DeclareMathOperator{\MHM}{MHM}
\DeclareMathOperator{\Map}{Map}
\DeclareMathOperator{\ev}{ev}
\DeclareMathOperator{\id}{id}
\newcommand{\R}{\mathbb{R}}
\newcommand{\Sc}{\mathcal{S}}
\newcommand{\Dol}{\mathrm{Dol}}
\newcommand{\Jac}{\mathrm{Jac}}
\newcommand{\Pic}{\mathrm{Pic}}
\newcommand{\tW}{\widetilde{W}}
\providecommand\sslash{/\!\!/}

\newcommand{\IH}{I\!H}
\newcommand{\IC}{IC}
\renewcommand{\AA}{\mathbb{A}}
\DeclareMathOperator{\Prym}{Prym}
\DeclareMathOperator{\codim}{codim}
\DeclareMathOperator{\spec}{spec}
\DeclareMathOperator{\nne}{end}

\calclayout

\title{Stable cohomology of universal character varieties}
\author{Ishan Banerjee}
\author{Faye Jackson}
\author{Anne Larsen}
\author{Sam Payne}
\author{Xiyan Zhong}
\date{\today}

\begin{document}

\begin{abstract}
    We study the \textit{universal $\PGL_n$-character variety} over $\M_{g}$, whose fiber over a point $[C]$ is the space of $\PGL_n$-local systems on the curve $C$, in the cases where these spaces are smooth. We use nonabelian Hodge theory and properties of Saito's mixed Hodge modules to show that the Leray--Serre spectral sequence for the projection to $\M_g$ degenerates at $E_2$.  As an application, we prove that the rational cohomology of these varieties stabilizes in the limit as $g \to \infty$ and compute the stable limit. We also deduce similar results for the universal $G$-character variety over $\M_{g,1}$ whose fiber over a punctured curve is the variety of $G$-local systems with fixed central monodromy around the puncture, for $G = \GL_n$ or $\SL_n$. The paper concludes with a computation of the ring structure on the stable cohomology and with an appendix by Anne Larsen and Mirko Mauri proving related results for the intersection cohomology of singular universal character varieties.
    \end{abstract}

\maketitle

\setcounter{tocdepth}{1}
\tableofcontents

\section{Introduction}

The degree $k$ cohomology group of the moduli space of curves of genus $g$ %
is independent of $g$ for $k$ less than or equal to $s(g) :=\frac{2g-2}{3}.$ These stable cohomology groups were computed by Madsen and Weiss \cite{Madsen-Weiss}, confirming a conjecture of Mumford \cite{Mumford83}.
Subsequently, Ellenberg, Venkatesh, and Westerland initiated an analogous study of cohomological stability for Hurwitz spaces, which parametrize branched $G$-covers of $\AA^1$ for a finite group $G$ or, equivalently, $G$-local systems on $\AA^1$ minus finitely many branch points, with prescribed local monodromy conditions  \cite{ellenberg-venkatesh-westerland}. %
More recently, Landesman and Levy computed the stable cohomology groups of these $G$-Hurwitz spaces \cite{landesmann-levy-stable-answer} as the number of branch points goes to $\infty$, and Putman \cite[Theorem C]{Putman-partial} proved a homological stability result for $G$-Hurwitz spaces of a genus $g$ curve with one boundary component.

Inspired by this recent progress, we investigate analogous  stability phenomena for the cohomology groups of universal character varieties that parametrize $G$-local systems on smooth curves of genus $g$, for certain algebraic groups $G$. We prove that these groups stabilize in the limit as $g \to \infty$ and compute the stable limits.

Our first such result is for $G = \PGL_n$.  We consider the universal character variety $\M_B^{\PGL}(\C_g,n) \to \M_g$ whose fiber over $[C]$ is the character variety $$[\Hom(\pi_1(C), \PGL_n(\CC))/\PGL_n(\CC)],$$ parametrizing $\PGL_n(\CC)$-local systems on $C$ via the monodromy representation. This character variety breaks into pieces indexed by $\ZZ/n\ZZ$, according to the obstruction to lifting a $\PGL_n(\CC)$-representation to $\SL_n(\CC)$, which is locally constant. We study the stable cohomology of the components $\M_B^{\PGL}(\C_g, n,d)$ corresponding to $d \in (\ZZ/n\ZZ)^*$, which are smooth. The appendix explains why this is the same as stable intersection cohomology for any $d \in \ZZ/n\ZZ$. See \S \ref{sec: Betti} for precise definitions and further details.

To describe the stable cohomology groups we use the following notation. Let 
\begin{equation} \label{eq:Vg}
    V_g := H^1(S_g;\QQ)
\end{equation}
be the first rational cohomology group of a closed orientable surface of genus $g$, viewed as a $\Mod_g$-module or, equivalently, a local system on $\M_g$.  Here, and throughout the main statements, we assume $g,n\geq2$. We consider $V_g$ as a graded object, supported in degree $1$, so $V_g[-i]$ is supported in degree $i + 1$. By convention, the grading on the cohomology of a graded local system is given by the sum of the cohomological degree and the degree of the local system, e.g. $H^k(\M_g;V_g[-2])$ is of degree $k+3$. We also define
\begin{equation}\label{eq: grading}
V_\circ := \QQ \oplus \QQ[-2], \quad \quad W_n := \QQ \oplus \QQ[-2] \oplus \cdots \oplus \QQ[2-2n], \quad \quad %
\quad  \widetilde{W}_n := W_n / \QQ. %
\end{equation}
Here, we think of the graded vector space $V_\circ$ as the even cohomology of a closed orientable surface $S_g$, supported in degrees 0 and 2. One may similarly think of $W_n$ as the rational cohomology of $\PP^{n-1}$.

Many of our results hold in a range of degrees, depending on $g$ and $n$, bounded above by 
\begin{equation}\label{eq:tautological-range}
R_{g,n}:=2(g-1)(n-1)+2.
\end{equation}
For $g,n\geq 2$, we have $s(g)\leq R_{g,n}$.

\begin{theorem} \label{thm:mainPGL}
For $d \in (\ZZ/n\ZZ)^*$ and  $k\le R_{g,n}$, the cohomology group $H^k(\M_B^{\PGL}(\C_g, n,d); \QQ)$ is isomorphic to the degree $k$ part of the graded vector space
\begin{equation}\label{eq: PGLresult}
H^*\big(\M_g; \textstyle \bigwedge^\bullet (V_g \otimes \tW_n)\big) \otimes \Sym^\bullet (V_\circ \otimes \tW_n).
\end{equation}
\end{theorem}

The following statement is a key technical step in the proof of Theorem~\ref{thm:mainPGL}, which is also of independent interest. Let $\pi \colon \M_B^{\PGL}(\C_{g},n,d) \to \M_g$ be the forgetful map taking a $\PGL_n$-local system on a smooth curve to the underlying curve.

\begin{theorem} \label{thm:E2PGL}
    Assume $d \in (\ZZ/n\ZZ)^*$. Then the Leray--Serre spectral sequence for the fibration $\pi \colon \M_B^{\PGL}(\C_{g},n,d) \to \M_g$ degenerates at $E_2$.
\end{theorem}

\noindent We prove this degeneration using the relative nonabelian Hodge correspondence \cite{Simpson-reps2}, purity of the cohomology of the moduli space of Higgs bundles \cite[Corollary~1.3.2]{Semiprojective}, and a splitting result for pure Hodge modules \cite[(4.5.4)]{Saito}. Theorem~\ref{thm:mainPGL} is deduced from the degeneration of the Leray--Serre spectral sequence for $\pi$ using previously known results on generators and relations of the cohomology of the fibers (\cite[Theorem~7]{Markman} and \cite[Lemma 4.1.12]{Hausel-RV}). %

Combining Theorem~\ref{thm:mainPGL} with results of Looijenga on the stabilization of cohomology of $\M_g$ with coefficients in symplectic local systems \cite{looijenga-symplectic-coeffs}, we deduce that the cohomology stabilizes with respect to $g$, and then also with respect to $n$, in degrees at most $s(g)$.

\begin{corollary} \label{cor:stab-PGL}
For $d \in (\ZZ/n\ZZ)^*$ and $k \leq s(g)$, the group $H^k(\M_B^{\PGL}(\C_g, n,d); \QQ)$ does not depend on $g$ and vanishes if $k$ is odd.  If, in addition, $k \le 2n-1$ then this group is independent of $n$. 
\end{corollary}

We furthermore compute the ring structure on the stable cohomology, as follows.   Regard $\QQ[\psi,c_2,\ldots,c_n]$ as a graded vector space with basis the monomials in the formal variables $\psi,c_2,\ldots,c_n$, where $\psi$ has degree $2$ and $c_j$ has degree $2j$. We define
\begin{equation}\label{eq:Kn}
K_n:=\big(\QQ[\psi,c_2,\ldots,c_n][2]\big)_{>0}.
\end{equation}
Then $K_n$ has a basis of formal symbols $\kappa_{a_0;\,a_2,\ldots,a_n}$,
indexed by nonnegative integer tuples $(a_0; a_2,\ldots,a_n)$ satisfying $2a_0+4a_2+\cdots+2na_n>2$, with
\[
\deg\big(\kappa_{a_0;\,a_2,\ldots,a_n}\big)=2a_0+4a_2+\cdots+2na_n-2.
\]

Consider the universal curve $p\colon\C\longrightarrow\M_B^{\PGL}(\C_g,n,d)$, and let $\psi_\C := c_1(\omega_p)$. Let $\mathcal P$ be the tautological principal $\PGL_n$-bundle on $\C$. For $2\leq j\leq n$, let $c_j(\mathcal P)$ denote its rational characteristic classes, with the convention of \S\ref{subsec: monodromy}.

\begin{theorem} \label{thm:stablePGL}
For $d \in (\ZZ/n\ZZ)^*$, the homomorphism of graded $\QQ$-algebras
\begin{equation}\label{eq: PGLring}
\Sym^\bullet (K_n) \longrightarrow H^*(\M_B^{\PGL}(\C_g, n,d); \QQ), \qquad \kappa_{a_0;\,a_2,\ldots,a_n} \longmapsto p_!\left(\psi_\C^{a_0}c_2(\mathcal P)^{a_2}\cdots c_n(\mathcal P)^{a_n}\right),
\end{equation}
is an isomorphism in degrees less than or equal to $s(g)$.
\end{theorem}

\noindent When no confusion seems possible, we also write $\kappa_{a_0;\,a_2,\ldots,a_n} \in H^*(\M_B^{\PGL}(\C_g, n,d); \QQ)$ for the image of the corresponding formal symbol under \eqref{eq: PGLring}. Theorem~\ref{thm:stablePGL} says that, in the stable range (degrees less than or equal to $s(g)$), the rational cohomology of $\M_B^{\PGL}(\C_g, n,d)$ is the polynomial algebra on these generalized $\kappa$-classes, directly analogous to the Madsen--Weiss description of the stable cohomology of $\M_g$. The dimensions of the stable cohomology groups of $\M_B^{\PGL}(\C_g, n, 1)$ for $n \leq 8$ in even degrees up to 20 are displayed in \Cref{fig:table-pgl}.

\medskip

We also consider universal $G$-character varieties over $\M_{g,1}$, for $G$ = $\GL_n$ or $\SL_n$.  The fiber of the universal $G$-character variety $M^{G}_B(\C_{g,1}, n, d)$ over a punctured curve $C$ is 
$$\{ \phi \in \Hom(\pi_1(C), G(\CC)): \phi(\gamma) = e^{2\pi i d/n} \Id\}/\PGL_n(\CC).$$
Here, $\gamma$ is a chosen loop around the puncture. Projectivization induces maps from these universal $G$-character varieties to the pullback $\M_B^{\PGL}(\C_g,n,d) \times_{\M_g} \M_{g,1}$, and the quotient and product descriptions below relate the cohomology of the three families.

\begin{theorem}\label{thm:sl-stability-answer}
For $d \in (\ZZ/n\ZZ)^*$ and $k\leq R_{g,n}$, the cohomology group  $H^k(M^{\SL}_B(\C_{g,1}, n, d);\QQ)$ is isomorphic to the degree $k$ part of the graded vector space
    \[
H^*(\M_{g,1}; \textstyle{\bigwedge^\bullet}(V_g \otimes \tW_n)) \otimes \Sym^\bullet(V_\circ \otimes \tW_n).
    \]
    Similarly, the cohomology group  $H^k(M^{\GL}_B(\C_{g,1}, n, d);\QQ)$ is isomorphic to the degree $k$ part of the graded vector space
    \[
H^*(\M_{g,1}; \textstyle{\bigwedge^\bullet}(V_g \otimes W_n)) \otimes \Sym^\bullet(V_\circ \otimes \tW_n).
    \]
    \end{theorem}

From Theorem~\ref{thm:sl-stability-answer}, we deduce stability with respect to $g$ for $k\leq s(g)$, and stability with respect to $n$ when, in addition, $k\leq2n-1$.
 
 \begin{corollary} \label{cor:stab-G}
 Assume $d \in (\ZZ/n\ZZ)^*$, and fix $G = \SL_n$ or $\GL_n$. For $k \leq s(g)$, the group $H^k(M^{G}_B(\C_{g,1},n,d);\QQ)$ does not depend on $g$ and vanishes if $k$ is odd. If, in addition, $k \le 2n-1$, then this group is independent of $n$. 
 \end{corollary}   
    
We compute the ring structure on the stable cohomology as follows.   We retain the notation $K_n$ and the generalized $\kappa$-classes introduced above.

For $G=\SL_n$ or $\GL_n$, let
\[
p_G \colon \C_G \longrightarrow M_B^G(\C_{g,1},n,d)
\]
be the universal curve and set $\psi_G:=c_1(\omega_{p_G})$.  The projectivization of the tautological local system extends across the marked point and gives a tautological principal $\PGL_n$-bundle $\mathcal P_G$ on $\C_G$.  In the $\GL_n$ case, the determinant local system also extends across the marked point; let $\mathcal D_{\GL_n}$ denote the resulting determinant line bundle on $\C_{\GL_n}$.

For $\SL_n$, the generalized $\kappa$-classes are exactly as for $\PGL_n$. For $\GL_n$, we consider additional generalized $\kappa$-classes obtained by multiplying by powers of $c_1(\mathcal{D}_{\GL_n})$ before pushing forward along $p_{\GL_n}$. We therefore define a formal graded vector space
\begin{equation}\label{eq:K'n}
K'_n:= \big(\QQ[\psi,c_1,c_2,\ldots,c_n][2]\big)_{>0}.
\end{equation}
Once again, $\psi$ is in degree $2$ and $c_j$ is in degree $2j$. Then $K'_n$ has a basis of formal symbols
\[
\kappa_{a_0,a_1;\,a_2,\ldots,a_n},
\]
indexed by nonnegative integer tuples $(a_0,a_1; a_2,\ldots,a_n)$ satisfying $2a_0+2a_1+4a_2+\cdots+2na_n>2$, with
\[
\deg\big(\kappa_{a_0,a_1;\,a_2,\ldots,a_n}\big)=2a_0+2a_1+4a_2+\cdots+2na_n-2.
\]
Because we are now working on curves with marked points, the universal curve $p_G$ has the universal marked section $$\sigma_G \colon M_B^G(\C_{g,1}, n, d) \to \C_G.$$ Pulling back products of $\psi_G$ and the Chern classes of $\mathcal P_G$ along $\sigma_G$ gives additional marked-point tautological classes on $M_B^G(\C_{g,1}, n, d)$.

\begin{theorem} \label{thm:sl-stable-ring}
For $d \in (\ZZ/n\ZZ)^*$, there are homomorphisms of graded $\QQ$-algebras
\begin{equation}\label{eq: SLring}
\Sym^\bullet(K_n) \otimes \QQ[\psi,c_2,\dots,c_n] \longrightarrow H^*(M^{\SL}_B(\C_{g,1},n,d);\QQ)
\end{equation}
and
\begin{equation}\label{eq: GLring}
\Sym^\bullet(K'_n) \otimes \QQ[\psi,c_2,\dots,c_n] \longrightarrow H^*(M^{\GL}_B(\C_{g,1},n,d);\QQ)
\end{equation}
that are isomorphisms in degrees $k \leq s(g)$. On the first tensor factor, the map in \eqref{eq: SLring} sends a formal $\kappa$-symbol to the Gysin push-forward of the corresponding monomial in $\psi_{\SL_n}$ and the classes $c_j(\mathcal P_{\SL_n})$, as in \eqref{eq: PGLring}. The map in \eqref{eq: GLring} sends $\kappa_{a_0,a_1;\,a_2,\ldots,a_n}$ to
\[
(p_{\GL_n})_!\left(\psi_{\GL_n}^{a_0}c_1(\mathcal D_{\GL_n})^{a_1}c_2(\mathcal P_{\GL_n})^{a_2}\cdots c_n(\mathcal P_{\GL_n})^{a_n}\right).
\]
On the second tensor factor, both maps are given by $\psi\mapsto\sigma_G^*(\psi_G)$ and $c_j\mapsto\sigma_G^*c_j(\mathcal P_G)$.
\end{theorem}

 \noindent There is no $c_1$ in the second tensor factor for $\GL_n$ because the determinant line bundle $\mathcal D_{\GL_n}$ is trivialized along the universal marked section.
 
 The dimensions of the stable cohomology groups of $M_B^{\SL}(\C_{g,1},n,1)$ and $M_B^{\GL}(\C_{g,1},n,1)$ for $n \leq 8$ in even degrees up to 20 are displayed in \Cref{fig:table-sl,fig:table-gl}.

\begin{remark}
Theorem~\ref{thm:stablePGL} gives a second proof of Corollary~\ref{cor:stab-PGL}, independent of the relative nonabelian Hodge correspondence and of the splitting theorem of \S\ref{sec:degeneration}. The proofs of the $\SL_n$- and $\GL_n$-ring formulas likewise avoid the splitting theorem. The $\SL_n$ proof uses the low-degree $\widetilde\Gamma$-invariance of Lemma~\ref{lemma: SL invariance}, whose proof appeals to the fiberwise nonabelian Hodge correspondence. The $\GL_n$ proof instead uses only the quotient description of Proposition~\ref{prop: relative B relations} and the identification of the $\widetilde\Gamma$-invariant part of the $\SL_n$ cohomology with the $\PGL_n$ cohomology; in particular, it does not use the nonabelian Hodge correspondence. None of these ring calculations identifies the Leray filtration or reproves the degeneration theorem, and they apply only in the smaller range $k\leq s(g)$.
\end{remark}

\begin{remark}
\Cref{thm:mainPGL,thm:sl-stability-answer} show that the rational cohomology groups of
\[
\M_B^{\PGL}(\C_g,n,d), \qquad
M_B^{\SL}(\C_{g,1},n,d), \qquad
M_B^{\GL}(\C_{g,1},n,d)
\]
are independent of $d\in(\ZZ/n\ZZ)^*$ in degrees $k\leq R_{g,n}$. In the stable range $k\leq s(g)$, \Cref{thm:stablePGL,thm:sl-stable-ring} give a simple explanation---the source algebras do not involve $d$.

For a fixed curve, the corresponding character varieties may be defined over $\QQ[e^{2\pi i/n}]$, and the varieties associated with the primitive classes $d\in(\ZZ/n\ZZ)^*$ are Galois conjugate. This gives isomorphisms of their $\ell$-adic étale cohomology, but it does not directly compare the universal total spaces, which are constructed as Borel quotients in \Cref{defn: relative total PGL,defn: relative gl}. Our proof in the larger range uses the relative nonabelian Hodge correspondence for $\PGL_n$. The statements for $\SL_n$ and $\GL_n$ are then deduced from the quotient and product descriptions of \Cref{prop: relative B relations}. The independence of $d$ in these cases has no evident direct geometric explanation on the Dolbeault side, although it belongs to the broader $\chi$-independence phenomena discussed in \cite[\S 0.3]{chi-indep}. The appendix proves that, for arbitrary $d\in\ZZ/n\ZZ$, the intersection cohomology is independent of $d$ in degrees $k<4(g-1)(n-1)-2$. This bound is optimal except when $n=2$.
\end{remark}

\begin{figure}[ht!]
\renewcommand{\arraystretch}{1.35}
\setlength{\tabcolsep}{5pt}
	\begin{tabular}{|c|c|c|c|c|c|c|c|c|c|c|c|}
		\hline
		degree $k$: & 0 & 2 & 4 & 6 & 8 & 10 & 12 & 14 & 16 & 18 & 20 \\
		\hline
		$\dim_\QQ H^{k}\big(\M_B^{\PGL}(\C_g,2,1); \QQ\big)$ & 1 & 2 & 5 & 11 & 23 & 45 & 87 & 160 & 290 & 512 & 889 \\
		$\dim_\QQ H^{k}\big(\M_B^{\PGL}(\C_g,3,1); \QQ\big)$ & 1 & 2 & 6 & 14 & 33 & 71 & 152 & 308 & 618 & 1201 & 2302 \\
		$\dim_\QQ H^{k}\big(\M_B^{\PGL}(\C_g,4,1); \QQ\big)$ & 1 & 2 & 6 & 15 & 36 & 81 & 180 & 381 & 795 & 1615 & 3227 \\
		$\dim_\QQ H^{k}\big(\M_B^{\PGL}(\C_g,5,1); \QQ\big)$ & 1 & 2 & 6 & 15 & 37 & 84 & 190 & 409 & 870 & 1800 & 3670 \\
		$\dim_\QQ H^{k}\big(\M_B^{\PGL}(\C_g,6,1); \QQ\big)$ & 1 & 2 & 6 & 15 & 37 & 85 & 193 & 419 & 898 & 1875 & 3857 \\
		$\dim_\QQ H^{k}\big(\M_B^{\PGL}(\C_g,7,1); \QQ\big)$ & 1 & 2 & 6 & 15 & 37 & 85 & 194 & 422 & 908 & 1903 & 3932 \\
		$\dim_\QQ H^{k}\big(\M_B^{\PGL}(\C_g,8,1); \QQ\big)$ & 1 & 2 & 6 & 15 & 37 & 85 & 194 & 423 & 911 & 1913 & 3960 \\
		\hline
	\end{tabular}
	\caption{The Betti numbers of $\M_B^{\PGL}(\C_g,n,1)$ for $2 \leq n \leq 8$, for even degrees $k \leq 20$ in the stable range $k \leq s(g)$.}
	\label{fig:table-pgl}
\end{figure}

\begin{figure}[ht!]
\renewcommand{\arraystretch}{1.35}
\setlength{\tabcolsep}{4.5pt}
	\begin{tabular}{|c|c|c|c|c|c|c|c|c|c|c|c|}
		\hline
		degree $k$: & 0 & 2 & 4 & 6 & 8 & 10 & 12 & 14 & 16 & 18 & 20 \\
		\hline
		$\dim_\QQ H^{k}\big(M_B^{\SL}(\C_{g,1},2,1); \QQ\big)$ & 1 & 3 & 9 & 22 & 51 & 109 & 225 & 443 & 849 & 1579 & 2874 \\
		$\dim_\QQ H^{k}\big(M_B^{\SL}(\C_{g,1},3,1); \QQ\big)$ & 1 & 3 & 10 & 27 & 69 & 163 & 372 & 809 & 1713 & 3518 & 7067 \\
		$\dim_\QQ H^{k}\big(M_B^{\SL}(\C_{g,1},4,1); \QQ\big)$ & 1 & 3 & 10 & 28 & 74 & 181 & 429 & 971 & 2140 & 4582 & 9599 \\
		$\dim_\QQ H^{k}\big(M_B^{\SL}(\C_{g,1},5,1); \QQ\big)$ & 1 & 3 & 10 & 28 & 75 & 186 & 447 & 1028 & 2305 & 5024 & 10724 \\
		$\dim_\QQ H^{k}\big(M_B^{\SL}(\C_{g,1},6,1); \QQ\big)$ & 1 & 3 & 10 & 28 & 75 & 187 & 452 & 1046 & 2362 & 5189 & 11169 \\
		$\dim_\QQ H^{k}\big(M_B^{\SL}(\C_{g,1},7,1); \QQ\big)$ & 1 & 3 & 10 & 28 & 75 & 187 & 453 & 1051 & 2380 & 5246 & 11334 \\
		$\dim_\QQ H^{k}\big(M_B^{\SL}(\C_{g,1},8,1); \QQ\big)$ & 1 & 3 & 10 & 28 & 75 & 187 & 453 & 1052 & 2385 & 5264 & 11391 \\
		\hline
	\end{tabular}
	\caption{The Betti numbers of $M_B^{\SL}(\C_{g,1},n,1)$ for $2 \leq n \leq 8$, for even degrees $k \leq 20$ in the stable range $k \leq s(g)$.}
	\label{fig:table-sl}
\end{figure}

\begin{figure}[ht!]
\renewcommand{\arraystretch}{1.35}
\setlength{\tabcolsep}{3.5pt}
	\begin{tabular}{|c|c|c|c|c|c|c|c|c|c|c|c|}
		\hline
		degree $k$: & 0 & 2 & 4 & 6 & 8 & 10 & 12 & 14 & 16 & 18 & 20 \\
		\hline
		$\dim_\QQ H^{k}\big(M_B^{\GL}(\C_{g,1},2,1); \QQ\big)$ & 1 & 5 & 22 & 79 & 260 & 785 & 2243 & 6089 & 15897 & 40068 & 98067 \\
		$\dim_\QQ H^{k}\big(M_B^{\GL}(\C_{g,1},3,1); \QQ\big)$ & 1 & 5 & 23 & 87 & 302 & 965 & 2922 & 8417 & 23335 & 62505 & 162651 \\
		$\dim_\QQ H^{k}\big(M_B^{\GL}(\C_{g,1},4,1); \QQ\big)$ & 1 & 5 & 23 & 88 & 310 & 1007 & 3106 & 9127 & 25841 & 70758 & 188378 \\
		$\dim_\QQ H^{k}\big(M_B^{\GL}(\C_{g,1},5,1); \QQ\big)$ & 1 & 5 & 23 & 88 & 311 & 1015 & 3148 & 9311 & 26555 & 73295 & 196809 \\
		$\dim_\QQ H^{k}\big(M_B^{\GL}(\C_{g,1},6,1); \QQ\big)$ & 1 & 5 & 23 & 88 & 311 & 1016 & 3156 & 9353 & 26739 & 74009 & 199350 \\
		$\dim_\QQ H^{k}\big(M_B^{\GL}(\C_{g,1},7,1); \QQ\big)$ & 1 & 5 & 23 & 88 & 311 & 1016 & 3157 & 9361 & 26781 & 74193 & 200064 \\
		$\dim_\QQ H^{k}\big(M_B^{\GL}(\C_{g,1},8,1); \QQ\big)$ & 1 & 5 & 23 & 88 & 311 & 1016 & 3157 & 9362 & 26789 & 74235 & 200248 \\
		\hline
	\end{tabular}
	\caption{The Betti numbers of $M_B^{\GL}(\C_{g,1},n,1)$ for $2 \leq n \leq 8$, for even degrees $k \leq 20$ in the stable range $k \leq s(g)$.}
	\label{fig:table-gl}
\end{figure}

\noindent\textbf{Acknowledgments.} This work grew out of discussions at the AIM workshop on ``Motives and mapping class groups.'' We thank AIM for the stimulating and supportive environment and our fellow workshop participants for many engaging conversations. We are grateful to Aaron Landesman for proposing this project and to Eduard Looijenga for invaluable input in its early stages. We thank Bhargav Bhatt and Davesh Maulik for helpful conversations related to this work, and Daniel Litt and Dan Petersen for comments on an earlier version of the paper. We are especially grateful to Oscar Randal-Williams for suggesting the homotopy-theoretic arguments used to compute the stable cohomology rings.

F.J. and A.L. were partially supported by NSF Graduate Research Fellowships under Grant Nos.~2140001 and~2141064, respectively. S.P. was partially supported by NSF grant DMS--2542134 and carried out portions of this work while visiting the Institute for Advanced Study, supported by the Charles Simonyi Endowment. X.Z. is grateful to the Max Planck Institute for Mathematics in Bonn for its hospitality.

\section{Universal character varieties} \label{sec: Betti}

In this section, we start by defining the universal $\PGL_n$-character variety over $\M_g$. We then define the universal $G$-character varieties over $\M_{g,1}$, for $G$ equal to $\GL_n$ or $\SL_n$. We also record the quotient and product descriptions relating these spaces. %

\medskip

Let $S_g$ be a closed oriented surface of genus $g \ge 2$, and let $\Mod_g := \pi_0(\Diff^+(S_g))$ be its mapping class group.  The group  $\PGL_n(\CC)$ acts by conjugation on the affine variety
\[
X_{g,n} := \Hom(\pi_1(S_g), \PGL_n(\CC)).
\]
Let $\mu_n \subset \CC^*$ be the group of $n$th roots of unity, which we identify with $\ZZ/n\ZZ$ via $d \mapsto e^{2\pi i d / n}$. By the short exact sequence %
$$1 \to \mu_n \to \SL_n(\CC) \to \PGL_n(\CC) \to 1,$$
the obstruction to lifting a representation $\phi \colon \pi_1(S_g) \to \PGL_n(\CC)$ to $\SL_n(\CC)$ is given by an element of $H^2(S_g; \mu_n) = \ZZ/n\ZZ$, and this obstruction is locally constant on $X_{g,n}$.  We write $X_{g,n,d}$ for the subvariety parametrizing homomorphisms with obstruction class $d \in \ZZ/n\ZZ$. Via the Riemann--Hilbert correspondence, $X_{g,n,d}$ parametrizes framed $\PGL_n$-bundles with flat connection that lift to $C^\infty$-vector bundles of degree $d$ mod $n$.

We now define the degree $d$ $\PGL_n$-character variety to be the quotient stack
\[
\M_{B}^{\PGL}(S_g,n,d) := [X_{g,n,d} / \PGL_n(\CC)],
\]
where $\PGL_n(\CC)$ acts by conjugation on the target. Likewise, we define the degree $d$ coarse $\PGL_n$-character variety to be the affine GIT quotient
\begin{equation} \label{eq:affineGIT-PGL}
M_B^{\PGL}(S_g,n,d) = X_{g,n,d} \sslash \PGL_n(\CC).
\end{equation}
In other words, $M_B^{\PGL}(S_g,n,d)$ is the affine variety whose coordinate ring is the ring of $\PGL_n(\CC)$-invariant regular functions on $X_{g,n,d}$. The subscript $B$ indicates that these character varieties appear on the Betti side of the nonabelian Hodge correspondence.

\begin{remark}
We write $\M$ for quotient stacks and $M$ for coarse GIT quotients. For the $\PGL_n$ character varieties, when the degree $d$ is coprime to $n$, the coarse-moduli morphism $\M \to M$ induces an isomorphism on rational cohomology, and we pass between the two using this isomorphism without further mention.
\end{remark}

Assume $d \in (\ZZ/n\ZZ)^*$. Then $X_{g,n,d}$ is smooth and irreducible, and $\PGL_n(\CC)$ acts with finite stabilizers \cite[Theorem~2.2.12 and Corollary~3.5.5]{Hausel-RV}.
In particular, $\M_{B}^{\PGL}(S_g,n,d)$ is a smooth orbifold, and projection to the coarse space $\M_{B}^{\PGL}(S_g,n,d) \to M_B^{\PGL}(S_g,n,d)$ induces an isomorphism %
\begin{equation}\label{eq:coarsespaceiso-PGL}
H^*(M_B^{\PGL}(S_g,n,d);\QQ) \cong H^*(\M_{B}^{\PGL}(S_g,n,d);\QQ).
\end{equation}

Let $\mc T_g$ denote the Teichm\"uller space of $S_g$, which parametrizes complex structures on $S_g$ up to the action of the identity component of $\Diff^+(S_g)$. It carries a natural action of $\Mod_g$, and the quotient is the moduli space of curves $\mc T_g / \Mod_g \cong \M_g$. The outer action of $\Mod_g$ on $\pi_1(S_g)$ also induces a $\Mod_g$-action on the character variety $\M_{B}^{\PGL}(S_g,n,d)$.

\begin{definition}\label{defn: relative total PGL}
	For $d$ coprime to $n$, the \emph{universal degree $d$ $\PGL$-character variety} is 
	\begin{equation*}
		\M_{B}^{\PGL}(\C_g, n, d) \coloneqq \mc T_g \times_{\Mod_g} \M_{B}^{\PGL}(S_g,n,d).
	\end{equation*}
\end{definition}

\noindent By construction, this universal character variety is a smooth and irreducible orbifold. It is a smooth fiber bundle over the moduli stack $\M_g$ with fiber $\M_{B}^{\PGL}(S_g, n, d)$.

\begin{remark}
Although $\M_g$ and each fiber are algebraic, the Borel construction above is only a topological orbifold; to apply techniques from algebraic geometry to study the rational cohomology of $\M_B^{\PGL}(\C_g, n, d)$, we pass to a homeomorphic algebraic moduli space of semistable Higgs bundles on the Dolbeault side of the nonabelian Hodge correspondence.
\end{remark}

We now consider degree $d$ $G$-character varieties for a surface with a puncture, for $G = \SL_n$ or $\GL_n$.
Let $S_{g,1}$ be an orientable surface of genus $g$ with one puncture, and let $\gamma$ denote a loop around the puncture. For $d \in \ZZ/n\ZZ$, we consider the affine variety
\begin{equation} \label{eq:YG}
Y^G_{g,n,d} := \{\phi \in \Hom(\pi_1(S_{g,1}),G(\CC)) : \phi(\gamma) = e^{2\pi i d/n}\Id_n\},
\end{equation}
on which $G(\CC)$ acts by conjugation.  If $d \in (\ZZ/n\ZZ)^*$ then  $Y^G_{g,n,d}$ is smooth and irreducible, and the $G(\CC)$-action factors through a free action of $\PGL_n(\CC)$ \cite[Corollary~2.2.7]{Hausel-RV}. In these cases, the degree $d$ $G$-character variety is the quotient
\begin{equation}
	M_B^G(S_{g,1}, n,d) \coloneqq Y^G_{g,n,d}/\PGL_n(\CC).%
    \label{eq:degree-d-single-curve}
\end{equation}
Unlike $\M_{B}^{\PGL}(S_g,n,d)$, which is an orbifold (or smooth Deligne--Mumford stack with nontrivial stabilizers), the space $M_B^G(S_{g,1}, n,d)$ is a manifold (or smooth algebraic variety).

Because the outer action of the mapping class group $\Mod_{g,1} \coloneqq \pi_0(\Diff^+(S_{g,1}))$ on $\pi_1(S_{g,1})$ preserves the conjugacy class of $\gamma$ \cite[Theorem 8.8]{Farb-Margalit}, $M_{B}^G(S_{g,1},n,d)$ inherits an action of the mapping class group $\Mod_{g,1}$ of $S_{g,1}$, by precomposition. Let $\mc T_{g,1}$ denote the Teichm\"uller space of $S_{g,1}$, which likewise carries a natural action of $\Mod_{g,1}$.

\begin{definition}%
\label{defn: relative gl}
	The universal degree $d$ $G$-character variety for $G = \GL_n$ or $\SL_n$ is
\begin{equation}
		M_{B}^G(\C_{g,1}, n, d) \coloneqq \mc T_{g,1} \times_{\Mod_{g,1}} M_{B}^G(S_{g,1},n,d).
	\end{equation}
\end{definition}

\begin{remark} \label{remark: point-pushing}
	When $d = 0$, the condition  $\phi(\gamma) = \Id_n$ from \eqref{eq:YG} is equivalent to the condition that $\phi$ factors through the fundamental group of the unpunctured surface $\pi_1(S_{g,1}) \to \pi_1(S_g).$ The outer action of $\Mod_{g,1}$ on $\pi_1(S_g)$ then factors through the outer action of $\Mod_g$, and $M_B^G(\C_{g,1},n,0)$
    is the pullback of an analogously defined universal degree 0 $G$-character variety over $\M_g$. Hence, one can study degree 0 relative $G$-character varieties for families of unpunctured curves, as in   \cite[\S 6 and \S9]{Simpson-reps2}.

However, if $d \in (\ZZ/n\ZZ)^*$ then $M_{B}^{G}(\C_{g,1},n,d)$ is not the pullback of a space over $\M_g$, because the kernel of the map $\pi_1(\M_{g,1}) \to \pi_1(\M_g)$ acts nontrivially on the fiber. This can be seen from the description of its outer action on $\pi_1(S_g)$  \cite[Fact 4.7]{Farb-Margalit}. 
\end{remark}

We define a relative degree $d$ $\PGL$-character variety over $\M_{g,1}$ (or indeed any family of curves) by pullback:
    \[
    \M_B^{\PGL}(\C_{g,1},n,d) := \M_B^{\PGL}(\C_g,n,d) \times_{\M_g} \M_{g,1}.
    \]
    This is analogous to the definition for $G = \GL_n$ or $\SL_n$, since the $\PGL$-character variety parametrizes homomorphisms that factor through $\pi_1(S_g)$.

We now discuss relations among the universal $G$-character varieties for varying groups $G$. Note that $\Gamma := \Hom(\pi_1(S_{g,1}), \mu_n)$ acts on $Y_{g,n,d}^G$ by multiplication, via the inclusion of $\mu_n$ in the center of $G$, and this action descends to an action on the character variety $M_B^G(S_{g,1},n,d)$. Passing to the universal setting, the relative group object $$ \widetilde\Gamma := \mc T_{g,1}  \times_{\Mod_{g,1}} \Gamma$$ acts on $M_B^G(\C_{g,1},n,d)$ over $\M_{g,1}$. %

\begin{prop} \label{prop: relative B relations}
Projectivization induces a map $q_{\SL}\colon M_B^{\SL}(\C_{g,1},n,d) \to \M_B^{\PGL}(\C_{g,1},n,d)$ 
which exhibits the target as the quotient stack
\[
\M_B^{\PGL}(\C_{g,1},n,d) \cong \big[M_B^{\SL}(\C_{g,1},n,d)/\widetilde\Gamma\big].
\]
In particular, on coarse spaces, $M_B^{\SL}(\C_{g,1},n,d)/\widetilde\Gamma \cong M_B^{\PGL}(\C_{g,1},n,d)$. Similarly, tensor product induces a quotient map $q_{\GL}\colon M_B^{\GL}(\C_{g,1},1,0) \times_{\M_{g,1}} M_B^{\SL}(\C_{g,1},n,d)
\to M_B^{\GL}(\C_{g,1},n,d)$,
and identifies its target with the quotient by the diagonal action of $\widetilde\Gamma$.
\end{prop}

\begin{proof}
    First, note that $\PGL_n = \SL_n/\mu_n$, where $\mu_n$ acts by its central embedding. Similarly,  $\GL_n = (\GL_1 \times \SL_n)/\mu_n$, where $\mu_n$ acts diagonally. Taking $\Hom(\pi_1(S_{g,1}), -)$ gives %
    \begin{equation*} %
        X_{g,n,d} = Y_{g,n,d}^{\SL}/ \Gamma 
    \quad \mbox{ and } \quad 
    Y^{\GL}_{g,n,d} = (Y_{g,1,0}^{\GL} \times Y_{g,n,d}^{\SL})/\Gamma.
    \end{equation*}
    In each case, the action of $\mu_n$ commutes with conjugation by $\PGL_n$ and with the action of $\Mod_{g,1}$ on $\pi_1(S_{g,1})$. Passing to the quotients and then to the Borel constructions gives the stated morphisms and quotient identifications. %
\end{proof}

\section{Degeneration of the Leray--Serre spectral sequence}\label{sec:degeneration}

In this section, we prove the degeneration of the Leray--Serre spectral sequence for the universal degree $d$ $\PGL$-character variety $\M_B^{\PGL}(\C_g, n,d) \to \M_g$. Similar results for universal degree $d$ $\GL$- and $\SL$-character varieties follow as corollaries. Our proof uses the relative nonabelian Hodge correspondence for $\PGL$. This relates the universal degree $d$ $\PGL_n$-character variety on the Betti side of the correspondence to the universal moduli space of semistable degree $d$ $\PGL_n$-Higgs bundles on the Dolbeault side. On the Dolbeault side, the projection to $\M_g$ shares many cohomological properties of a smooth and proper morphism of algebraic varieties, despite not being proper. 
Using these properties, we prove the following strong form of Theorem~\ref{thm:E2PGL}.%

\begin{theorem} \label{thm: splitting}
Let $d \in (\ZZ/n\ZZ)^*$ and let $\pi \colon \M_B^{\PGL}(\C_g, n,d) \to \M_g$ be the projection from the universal degree $d$ $\PGL$-character variety. Then $R\pi_* \QQ \cong \oplus_i R^i\pi_* \QQ[-i]$ in the bounded derived category of constructible $\QQ$-sheaves on $\M_g$. Hence, the Leray--Serre spectral sequence
    $$E_2^{pq} = H^p(\M_g, R^q\pi_* \QQ) \implies H^{p+q}(\M_B^{\PGL}(\C_g, n,d); \QQ)$$
    degenerates at $E_2$.
\end{theorem}

\begin{remark}
The first statement in Theorem~\ref{thm: splitting} is a universal form of the degeneration of the Leray--Serre spectral sequence in the following sense: an object $A$ in a bounded derived category $D^b(\mc A)$ splits as the direct sum of its cohomology sheaves if and only if for \textit{any} cohomological functor $F$ from $D^b(\mc A)$ to an abelian category, the corresponding spectral sequence $E_2^{pq} = F(\mc H^q(A)[p]) \implies F(A[p+q])$ degenerates at $E_2$ \cite[Proposition~1.2]{DeligneE2}.
\end{remark}


\begin{remark}
An earlier version of this paper contained a second, longer proof of the above theorem using only classical Hodge theory and some manipulations in the derived category of constructible sheaves \cite[Section~3.5]{BJLPZv2}. One may also give a third proof of Theorem~\ref{thm: splitting}, similar in spirit but different in technical detail, via spreading out over a finitely generated $\ZZ$-algebra, reducing modulo large finite primes, and applying the theory of weights in the $\ell$-adic cohomology of local systems on varieties over finite fields, along the lines of \cite[Theorem~5.4.5]{BBD}. Like the proof presented here, this argument relies on the theory of weights in the cohomology of local systems on algebraic varieties \cite{DelignePoids}. However, to present this third proof in detail would take us too far from the setting of this paper, which is rooted in geometry over the complex numbers. The second proof is also omitted from the present version, in the interest of efficiency.
\end{remark}

In \S \ref{subsec: PGL B to GL/SL}, we deduce the following corollaries:

\begin{corollary} \label{cor: splitting, GL}
    Let $d \in (\ZZ/n\ZZ)^*$, and let  $\pi\colon M_B^{\GL}(\C_{g,1},n,d) \to \M_{g,1}$ be the projection map from the universal degree $d$ $\GL_n$-character variety.  %
    Then $R\pi_* {\QQ} \cong \oplus R^i \pi_* {\QQ}[-i]$, so the Leray--Serre spectral sequence corresponding to $\pi$ degenerates at $E_2$.
\end{corollary}

For $n\geq 2$, let $p_n$ denote the smallest prime divisor of $n$, and set
\begin{equation}\label{eq:rho-range}
\rho_{g,n}:=2n^2(1-p_n^{-1})(g-1).
\end{equation}
Since $p_n\geq 2$, for $g,n\geq 2$ we have
\begin{equation}\label{eq:rho-bound}
\rho_{g,n}\geq n^2(g-1)\geq R_{g,n}.
\end{equation}

\begin{corollary} \label{cor: splitting, SL}
    Let $d \in (\ZZ/n\ZZ)^*$, and let $\pi\colon M_B^{\SL}(\C_{g,1},n,d) \to \M_{g,1}$ be the projection from the universal degree $d$ $\SL$-character variety. Then the truncation is $$\tau_{\le \rho_{g,n}} R\pi_* {\QQ} \cong \bigoplus_{i=0}^{\rho_{g,n}} R^i\pi_* {\QQ}[-i].$$
    In particular, the Leray--Serre spectral sequence for $\pi$ degenerates at $E_2$ in degrees $\le \rho_{g,n}$.
\end{corollary}

\subsection{The  moduli spaces}\label{subsec: relative dolbeault}
We now introduce the moduli spaces of Higgs bundles that appear on the Dolbeault side of the nonabelian Hodge correspondence. %

Recall that a \textit{Higgs bundle} on a curve $C$ is a pair $(E, \theta)$, where $E$ is a vector bundle on $C$ and $\theta\colon E \to E \otimes \omega_C$ a twisted endomorphism. The slope of a vector bundle $E$ is 
\[
\mu(E) := \mathrm{degree} (E) / \mathrm{rank} (E).
\]
A Higgs bundle $(E, \theta)$ is \emph{stable} or  \emph{semistable} if, for every $\theta$-invariant subbundle $F \subset E$, the slopes satisfy $\mu(F) < \mu(E)$ or $\mu(F) \leq \mu(E)$, respectively. When the degree and rank are coprime, a Higgs bundle is stable if and only if it is semistable.   Let $M^{\GL}_{\Dol}(C,n,d)$ denote the coarse moduli space of semistable Higgs bundles of degree $d$ and rank $n$ on $C$, as constructed by Simpson \cite{Simpson-reps2}.

For moduli of $\SL$- and $\PGL$-Higgs bundles, we follow the definitions from \cite[\S 2]{Hausel-Thaddeus-mirror}. 
Fixing an auxiliary degree $d$ line bundle $L$ on $C$, the moduli space of $\SL$-Higgs bundles $M^{\SL}_{\Dol}(C,n,L)$ is the subspace of $M^{\GL}_{\Dol}(C,n,d)$ parametrizing Higgs bundles $(E, \theta)$ such that $\det E \cong L$ and $\Tr \theta = 0$, when  $\theta$ is viewed as an $\omega_C$-twisted endomorphism of $E$. 

\begin{definition}
The moduli space of $\PGL$-Higgs bundles is the geometric quotient 
\begin{equation} \label{eq:MPGL-Dol}
M^{\PGL}_{\Dol}(C,n,d) 
\cong M^{\SL}_{\Dol}(C,n,L)/\Jac(C)[n].
\end{equation}
\end{definition}

\noindent Here, the $n$-torsion subgroup $\Jac(C)[n]$ acts on $M^{\SL}_{\Dol}(C,n,L)$ by $M \cdot (E, \theta) := (M \otimes E, \theta).$ For a different choice of degree $d$ line bundle $L'$, there is a noncanonical isomorphism $M^{\SL}_{\Dol}(C,n,L) \cong M^{\SL}_{\Dol}(C,n,L')$ given by tensoring with an $n^{\mathrm{th}}$ root of $L^{-1} \otimes L' \in \Jac(C)$. It follows that $M^{\PGL}_{\Dol}(C,n,d)$ is well-defined by \eqref{eq:MPGL-Dol} and independent of $L$.

These Higgs moduli spaces are related to the character varieties on the Betti side of the nonabelian Hodge correspondence, as follows \cite[Theorem 4.3]{Hausel-handbook}.

\begin{NAHC}
    Let $C$ be a smooth projective curve of genus $g$. There is a homeomorphism 
    $M^{\PGL}_{\Dol}(C,n,d) \cong M^{\PGL}_B(S_g,n,d)$, and also 
    $$M^{\GL}_{\Dol}(C,n,d) \cong M^{\GL}_B(S_{g,1},n,d), \quad \mbox{ and } \quad M^{\SL}_{\Dol}(C,n,L) \cong M^{\SL}_B(S_{g,1},n,\deg L).$$
\end{NAHC}

\subsection{Smoothness and purity for \texorpdfstring{$d$}{d} coprime to \texorpdfstring{$n$}{n}}

We now switch to the moduli stack.

\begin{definition}
If $d \in (\ZZ/n\ZZ)^*$ then
$\M^{\PGL}_{\Dol}(C,n,d) := [M^{\SL}_{\Dol}(C,n,L)/\Jac(C)[n]].$
\end{definition}

Equivalently, there is an isomorphism
\begin{equation} \label{eq:PGL-Higgs from GL}
\M^{\PGL}_{\Dol}(C,n,d) \cong [M_{\Dol}^{\GL}(C,n,d)/M^{\GL}_{\Dol}(C,1,0)].
\end{equation}
Here,
\[
M^{\GL}_{\Dol}(C,1,0) \cong \Jac(C) \times H^0(C,\omega_C)
\]
acts by
\[
(M,s)\cdot(E,\theta)=(M\otimes E,\theta+s\Id_E).
\]
Thus $\M^{\PGL}_{\Dol}(C,n,d)$ is a Deligne--Mumford stack with coarse space $M^{\PGL}_{\Dol}(C,n,d)$. Unlike the coarse space, the moduli stack is smooth: by \eqref{eq:PGL-Higgs from GL}, it is the quotient of the smooth variety $M^{\GL}_{\Dol}(C,n,d)$ \cite[Proposition~7.4]{Nitsure} by the smooth group $M^{\GL}_{\Dol}(C,1,0)$. Since $M^{\SL}_{\Dol}(C,n,L)\to\M^{\PGL}_{\Dol}(C,n,d)$ is a $\Jac(C)[n]$-torsor, $M^{\SL}_{\Dol}(C,n,L)$ is also smooth. 

The Higgs moduli space $\M^{\PGL}_{\Dol}(C,n,d)$ is not proper, but it is proper over an affine space and shares some of the cohomological properties of a smooth and proper Deligne--Mumford stack, as we now recall.

The Hitchin morphism $M^{\SL}_{\Dol}(C,n,L) \to \bigoplus_{i=2}^n H^0(C, \omega_C^{\otimes i})$ takes a Higgs bundle $(E,\theta)$ to the coefficients of the characteristic polynomial of $\theta$, viewed as a twisted endomorphism of $E$. The Hitchin morphism is proper \cite[Theorem 6.11]{Simpson-reps2}. It is also equivariant with respect to the following $\mathbb{G}_m$-actions. The action on the source is $t \cdot (E, \theta) := (E, t\theta)$. On the target, for $x \in H^0(C, \omega_C^{\otimes i})$, the action is $t \cdot x := t^i x$. 
Because the $\mathbb{G}_m$-action contracts the base to the point 0 and the Hitchin morphism is proper, by relating the cohomology of $M^{\SL}_{\Dol}(C,n,L)$ to that of the proper $\mathbb{G}_m$-fixed locus $M_{\Dol}^{\SL}(C,n,L)^{\mathbb{G}_m}$, we have the following:

\begin{prop} \label{prop: purity}
\cite[Corollary~1.3.2]{Semiprojective}
If $d = \deg L \in (\ZZ/n\ZZ)^*$, then both of
$$H^k(M^{\SL}_{\Dol}(C,n,L); \QQ) \supset H^k(M^{\SL}_{\Dol}(C,n,L); \QQ)^{\Jac(C)[n]} = H^k(\M^{\PGL}_{\Dol}(C,n,d); \QQ)$$
are pure of weight $k$.
\end{prop}

\begin{remark}
The character variety $M_B^{\PGL}(C,n,d)$ is homeomorphic to $M^{\PGL}_{\Dol}(C,n,d)$, by the nonabelian Hodge correspondence, but this homeomorphism is not algebraic and the induced isomorphism on cohomology does not preserve weights.  Indeed, the weight filtration on $H^k(M_B^{\PGL}(C,n,d);\QQ)$ is typically far from pure, and is the subject of the $P=W$ conjecture proposed in \cite{PW}. The $P=W$ conjecture for $\PGL_n$ follows from the corresponding conjecture for $\GL_n$, which is proved in \cite{MS, HMMS}.
\end{remark}

\subsection{Universal Higgs moduli spaces}
The definition of $M^{\GL}_{\Dol}(C,n,d)$ readily generalizes to the relative setting for a family of curves over a base. The quotient presentation \eqref{eq:PGL-Higgs from GL} likewise generalizes to families. Thus, from the universal curve over $\M_g$, we obtain a universal $\PGL_n$-Higgs moduli space $M^{\PGL}_{\Dol}(\mc C_g,n,d) \to \M_g.$ %
By \cite[Theorem~9.11 and Lemma~9.14]{Simpson-reps2}, there is a homeomorphism from the relative $\PGL$-character variety to the relative $\PGL$-Higgs moduli space for an arbitrary family of curves. Simpson's homeomorphism preserves the $C^\infty$-isomorphism class of the underlying $\PGL$-bundles, and hence maps degree $d$ character varieties to degree $d$ Higgs moduli spaces.

\begin{RNAHC} %
There is a homeomorphism  over $\M_g$ $$M_B^{\PGL}(\C_g, n, d) \cong M_\Dol^{\PGL}(\C_g,n,d).$$ %
\end{RNAHC}

\begin{corollary} \label{cor: switch to Dol}
    Let $\pi_B \colon \M^{\PGL}_B(\C_g,n,d) \to \M_g$ be the universal degree $d$ $\PGL$-character variety. Then the pushforward $R\pi_{B*} \QQ$ %
    splits as a sum of its cohomology sheaves if and only if the same is true for 
	$
		\pi_{\Dol}\colon \M^{\PGL}_{\Dol}(\C_g,n,d) \to \M_g. %
	$
\end{corollary}

\begin{proof}
    We note that the pushforward of the constant sheaf $\QQ$ from a Deligne--Mumford stack to its coarse space is the constant sheaf $\QQ$. Factoring the projection
    $$\pi_B\colon \M_B^{\PGL}(\C_g,n,d) \xrightarrow{f_B} M_B^{\PGL}(\C_g,n,d) \xrightarrow{g_B} \M_g,$$
    we see $R\pi_{B*} \QQ \cong Rg_{B*}\QQ$.
Doing the same on the Dolbeault side shows $R\pi_{\Dol *} \QQ \cong Rg_{\Dol *} \QQ.$ The corollary follows, since $M_B^{\PGL}(\C_g, n, d) \cong M_\Dol^{\PGL}(\C_g,n,d)$ over $\M_g$.
\end{proof}

\subsection{Proof of Theorem \ref{thm: splitting}} \label{subsec: Saito}

To handle constructible sheaves equipped with Hodge structures, we use Saito's theory of mixed Hodge modules  \cite{Saito}. The key features are as follows. There is a category of mixed Hodge modules $\MHM(X)$ on each complex variety $X$, with a standard six-functor formalism. There is a faithful, exact functor
\begin{equation} \label{eq:rat}
\mathrm{rat}\colon D^b \MHM(X) \to D^b_c(\QQ_X)
\end{equation}
taking a mixed Hodge module to its underlying rational structure, %
which is compatible with the six functors \cite[Theorem 0.1]{Saito}. Moreover, $\MHM(\mathrm{pt})$ is the category of polarizable $\QQ$-mixed Hodge structures \cite[(4.2.12)]{Saito}, and pushing forward the trivial Hodge module $\QQ^H_X$ to a point we recover Deligne's mixed Hodge structure on $H^*(X;\QQ)$  \cite[p.~328]{Saito}.

The needed input from the theory of mixed Hodge modules, as extended to complex algebraic stacks in \cite{Tubach}, is the following lemma:

\begin{lemma}\label{lemma: MHM input}
    If $\pi \colon U \to S$ is a representable smooth morphism of smooth complex finite type Deligne--Mumford stacks such that the cohomology of each fiber $H^k(U_s; \QQ)$ has a pure Hodge structure of weight $k$ for each degree $k$ and $s \in S$, and such that each cohomology sheaf $R^i\pi_* \QQ_U$ is a local system on $S$, then
    $R\pi_* \QQ_U \cong \oplus_i R^i\pi_* \QQ_U[-i]$. In particular, the corresponding Leray--Serre spectral sequence degenerates at $E_2$.
\end{lemma}

\begin{proof}
    We start with the case of schemes.
    Our goal is to show that the (derived) pushforward $\pi_* \QQ^H_U \in D^b\MHM(S)$ is \emph{pure of weight 0} \cite[\S 4.5]{Saito}, i.e., that $\mathrm{Gr}^W_i \mc H^j(\pi_* \QQ^H_U) = 0$ for $i \ne j$. This condition guarantees that $\pi_* \QQ^H_U \cong \oplus_i \mc H^i(\pi_* \QQ^H_U)[-i] \in D^b \MHM(S)$ \cite[(4.5.4)]{Saito}, and then applying the rat functor (\ref{eq:rat}), one recovers the same identity on the level of constructible sheaves.

    To prove the purity condition, it suffices to check that for every $s \in S$ and $k \in \ZZ$ we have that $\mc H^k i_s^* \pi_* \QQ^H_U$ is of weight $\le k$ and $\mc H^k i_s^! \pi_* \QQ_U^H$ is of weight $\ge k$ \cite[(4.6.1)]{Saito}. We begin by establishing a few identities. First, recall the base-change identity \cite[(4.4.3)]{Saito}:
    \begin{equation} \label{eq: MHM base change}
        i_s^! \pi_* \QQ^H_U \cong \pi_{U_s *} i^!_{U_s} \QQ_U^H.
    \end{equation}
    Next, since $U$ and $\pi$ (therefore also $U_s$) are smooth, we have \cite[(4.4.2)]{Saito}
    \begin{equation} \label{eq: MHM smooth pullback}
      i_{U_s}^! \QQ^H_U \cong i^!_{U_s} p_U^! \QQ_{pt}^H(-d_U)[-2d_U] = p_{U_s}^! \QQ^H_{pt}(-d_U)[-2d_U] \cong \QQ^H_{U_s}(-d)[-2d], 
    \end{equation}
    where we write $d_U := \dim U$ and $d := \dim S$ and $p_U, p_{U_s}$ for the projection to a point.
    Finally, we consider the natural morphism $i_s^!(\QQ^H_S) \otimes i_s^*(-) \to i_s^!(-)$ obtained from adjunction and the projection formula; cf.~\cite[Proposition~3.1.11]{Kashiwara-Schapira} for the corresponding morphism of constructible complexes. Applying it to $\pi_* \QQ^H_U$ and taking the rat functor gives
    $$i_s^!(\QQ_S) \otimes i_s^*(R\pi_* \QQ_U) \longrightarrow i_s^!(R\pi_* \QQ_U).$$
    Since each $R^i\pi_* \QQ_U$ is a local system on $S$, the inclusion $i_s$ is non-characteristic for $R\pi_*\QQ_U$, so this morphism is an isomorphism by \cite[Proposition~5.4.13]{Kashiwara-Schapira}. Since the rat functor is faithful and exact, the original morphism in $D^b\MHM(S)$ is an isomorphism as well.
    By smoothness of $S$ we have $i_s^!(\QQ^H_S) = i_s^!(p_S^! \QQ^H_{pt}(-d)[-2d]) = \QQ_s^H(-d)[-2d]$ \cite[(4.4.2)]{Saito}. Putting these together, we obtain an isomorphism
    \begin{equation} \label{eq: MHM pullback}
        i_s^*(\pi_* \QQ_U^H)(-d)[-2d] \cong i_s^!(\pi_* \QQ_U^H).
    \end{equation}
    
    With these preliminaries, for the first condition we have
    $$i_s^* \pi_* \QQ^H_U \overset{(\ref{eq: MHM pullback})}{\cong} i_s^! \pi_* \QQ^H_U(d)[2d] \overset{(\ref{eq: MHM base change})}{\cong} \pi_{U_s *} i^!_{U_s} \QQ^H_U(d)[2d] \overset{(\ref{eq: MHM smooth pullback})}{\cong} \pi_{U_s *} \QQ^H_{U_s} \cong H^*(U_s; \QQ)$$
    which is of weight 0 by assumption.
    Similarly, for the second condition
    $$i_s^! \pi_* \QQ^H_U \overset{(\ref{eq: MHM base change})}{\cong} \pi_{U_s *} i^!_{U_s} \QQ^H_U \overset{(\ref{eq: MHM smooth pullback})}{\cong} \pi_{U_s *} \QQ^H_{U_s}(-d)[-2d] \cong H^{* - 2d}(U_s; \QQ)(-d),$$
    which is again of weight 0 by assumption.

    Finally, returning to the case of stacks, by \cite[Corollary 3.24]{Tubach} it again suffices to show that $\pi_* \QQ^H_U \in D^b \MHM(S)$ is pure of weight 0 \cite[Definition 3.19]{Tubach}. Taking a smooth cover by a variety $s:\overline{S} \to S$, as in \cite[Proof of Proposition 3.21]{Tubach} it suffices to check that $s^* \pi_* \QQ^H_U$ is pure of weight 0. By smooth base change \cite[Theorem 3.1(4)]{Tubach} we have $s^* \pi_* \QQ^H_U \cong \bar{\pi}_* \QQ^H_{U \times_S \overline{S}}$, for $\bar{\pi}: U \times_S \overline{S} \to \overline{S}$, and as $\pi$ was assumed representable, we now apply the result for schemes.
\end{proof}

\begin{proof}[Proof of Theorem~\ref{thm: splitting}]
	By Corollary~\ref{cor: switch to Dol} it suffices to prove the statement on the Dolbeault side, for $\pi \colon \M_{\Dol}^{\PGL}(\C_g,n,d) \to \M_g$. Let $\widetilde{M}_{\Dol}^{\PGL}(\C_g, n,d) \subset M_{\Dol}^{\GL}(\C_g,n,d)$ be the closed substack parametrizing Higgs bundles $(E, \theta)$ on $C$ such that $(\det E)^{\otimes (2g-2)} \cong \omega_C^{\otimes d}$, so that fiberwise $\widetilde{M}_{\Dol}^{\PGL}(C,n,d)$ is the disjoint union of $(2g-2)^{2g}$ copies of $M_{\Dol}^{\SL}(C,n,L)$ ranging over $L$ with $L^{\otimes (2g-2)} = \omega_C^{\otimes d}$. Note that $\widetilde{M}_{\Dol}^{\PGL}(\C_g,n,d) \to \M_{\Dol}^{\PGL}(\C_g,n,d)$ is a quotient by the finite group object
    \[ 
    \Delta:=\Pic_{\C_g/\M_g}^0[(2g-2)n]
    \]
    over $\M_g$.

    Consider the projection $\widetilde\pi: \widetilde{M}_{\Dol}^{\PGL}(\C_g,n,d) \to \M_g$. Recall that $\widetilde\pi$ is a smooth morphism, by \cite[Proposition 7.4]{Nitsure}, and its fibers have pure cohomology, by   Proposition~\ref{prop: purity}. We note that $\widetilde\pi$ is topologically locally trivial by the Universal Nonabelian Hodge Correspondence for $\PGL$ combined with the Fiberwise Nonabelian Hodge Correspondence for $\SL$ (see the proof of \cite[Proposition 3.6]{SLmonodromy} for more discussion). Thus in particular the cohomology sheaves $R^i \widetilde\pi_* \QQ$ are locally free. It follows from Lemma~\ref{lemma: MHM input} that
    \[
    R\widetilde\pi_*\QQ \cong \bigoplus_i R^i\widetilde\pi_*\QQ[-i].
    \]
    Since $\Delta$ is finite, taking invariants with rational coefficients is exact. Thus
    \[
    R\pi_*\QQ
    \cong (R\widetilde\pi_*\QQ)^\Delta
    \cong \bigoplus_i (R^i\widetilde\pi_*\QQ)^\Delta[-i]
    \cong \bigoplus_i R^i\pi_*\QQ[-i]. \qedhere
    \]
\end{proof}

\subsection{Proof of Corollaries \ref{cor: splitting, GL} and \ref{cor: splitting, SL}: the case of \texorpdfstring{$\GL$}{GL} and \texorpdfstring{$\SL$}{SL}} \label{subsec: PGL B to GL/SL}
Here, we deduce the splitting results for the $\GL$ and $\SL$ relative character varieties from the result for $\PGL$ (Theorem~\ref{thm: splitting}), using the morphisms described in Proposition~\ref{prop: relative B relations}. For $\SL$, the key point is the low-degree invariance recorded in the following lemma. We deduce the result for $\GL$ from the $\SL$ result and a relative Künneth formula (Lemma~\ref{lemma: smooth relative Künneth}).

\begin{lemma}\label{lemma: SL invariance}
Let $d\in(\ZZ/n\ZZ)^*$, and let $q\colon M_B^{\SL}(\C_{g,1},n,d)\longrightarrow\M_B^{\PGL}(\C_{g,1},n,d)$  be the quotient by the relative group object $\widetilde\Gamma$ of Proposition~\ref{prop: relative B relations}. Write $\pi_{\SL}$ and $\pi_{\PGL}$ for the projections to $\M_{g,1}$. Then $q$ induces an isomorphism
\begin{equation}\label{eq:RpiGamma}
R\pi_{\PGL *}\QQ \cong (R\pi_{\SL *}\QQ)^{\widetilde\Gamma}.
\end{equation}
Moreover, for $i \le \rho_{g,n}$, the action of $\widetilde\Gamma$ on $R^i\pi_{\SL *}\QQ$ is trivial, and hence
$R^i\pi_{\PGL *}\QQ \cong R^i\pi_{\SL *}\QQ$. In particular, for $g,n\geq 2$, these conclusions hold in degrees at most $R_{g,n}$.
\end{lemma}

\begin{proof}
The quotient description and exactness of invariants with rational coefficients give \eqref{eq:RpiGamma}. Fix a pointed curve $(C,p)$. The fiber of $\widetilde\Gamma$ over $(C,p)$ is
\[
\Gamma_C := \Hom(\pi_1(C \smallsetminus \{p\}),\mu_n) \cong \Jac(C)[n],
\]
and taking cohomology of the fiber in \eqref{eq:RpiGamma} gives
\[
H^i(\M_B^{\PGL}(S_g,n,d);\QQ) \cong H^i(M_B^{\SL}(S_{g,1},n,d);\QQ)^{\Gamma_C}.
\]
Choose a line bundle $L$ on $C$ whose degree represents $d$. By \cite[Theorem~1.4]{Hausel-Pauly}, $\Jac(C)[n]$ acts trivially on $H^i(M_{\Dol}^{\SL}(C,n,L);\QQ)$ for $i \le \rho_{g,n}$. By the transfer for the finite quotient $M_{\Dol}^{\PGL}(C,n,d) \cong M_{\Dol}^{\SL}(C,n,L)/\Jac(C)[n]$, the cohomology of $M_{\Dol}^{\PGL}(C,n,d)$ is the $\Jac(C)[n]$-invariant part of the cohomology of $M_{\Dol}^{\SL}(C,n,L)$. It follows that
\[
\dim H^i(M_{\Dol}^{\PGL}(C,n,d);\QQ) = \dim H^i(M_{\Dol}^{\SL}(C,n,L);\QQ)
\]
in this range. The nonabelian Hodge correspondence gives the same equality for the corresponding Betti moduli spaces. Hence the invariant subspace $H^i(M_B^{\SL}(S_{g,1},n,d);\QQ)^{\Gamma_C}$ has the same dimension as $H^i(M_B^{\SL}(S_{g,1},n,d);\QQ)$, so $\Gamma_C$ acts trivially. This holds for every pointed curve, and therefore $\widetilde\Gamma$ acts trivially on the local system $R^i\pi_{\SL *}\QQ$. The final assertions follow from \eqref{eq:RpiGamma} and \eqref{eq:rho-bound}.
\end{proof}

\begin{proof}[Proof of Corollary \ref{cor: splitting, SL}]
First, recall that the splitting in Theorem~\ref{thm: splitting} is stated for the universal $\PGL_n$-character variety over $\M_g$. By smooth base change \cite[Theorem 3.2.13(ii)]{Dimca}, the analogous splitting holds for any family of curves induced by a smooth morphism to $\M_g$. In particular, it holds for the pullback to $\M_{g,1}$.

By Lemma~\ref{lemma: SL invariance}, pullback along $q$ induces an isomorphism on the cohomology sheaves in degrees at most $\rho_{g,n}$, and hence an isomorphism after truncation. Applying the base-changed $\PGL$ splitting, and then Lemma~\ref{lemma: SL invariance} once more, gives
\[
\tau_{\le \rho_{g,n}}R\pi_{\SL *}\QQ \cong \tau_{\le \rho_{g,n}}R\pi_{\PGL *}\QQ \cong \bigoplus_{i=0}^{\rho_{g,n}}R^i\pi_{\PGL *}\QQ[-i] \cong \bigoplus_{i=0}^{\rho_{g,n}}R^i\pi_{\SL *}\QQ[-i],
\] as required.
\end{proof}

In the proof of Corollary~\ref{cor: splitting, GL}, we use the following lemma.

\begin{lemma}\label{lemma: smooth relative Künneth}
    Let $\pi_i\colon X_i \to S$, $i=1,2$, be submersions of complex manifolds, and let $\pi\colon X_1 \times_S X_2 \to S$ be the induced projection. Suppose that each $R^j\pi_{2*} \QQ$ is locally trivial. Then $$R\pi_*\QQ \cong R\pi_{1*} \QQ \otimes R\pi_{2*} \QQ$$ in $D^b_c(\QQ_S)$. The same holds if the base $S$ is a smooth complex Deligne--Mumford stack.
\end{lemma}

\begin{proof}
    All of the functors in this proof are derived; for simplicity we omit this from the notation. Also, in this proof (only) we use $\dim$ to denote $\dim_{\R}$. We start with the case that $S$ is a manifold. The strategy is to apply Verdier duality to the proper relative K\"unneth formula \cite[Exercise II.18(i)]{Kashiwara-Schapira}
    \begin{equation} \label{eq:rel-Kuenneth}
    \pi_!(\QQ) = \pi_{1!} \QQ \otimes \pi_{2!} \QQ.
    \end{equation}
    On the left side we have
    $$D\pi_!\QQ = \pi_* D \QQ = \pi_*\QQ[\dim X_1 \times_S X_2]$$
    \cite[Theorem 3.2.17(iii) and Proposition 3.3.7]{Dimca}. %
    On the right side we have
    $$D(\pi_{1!} \QQ \otimes \pi_{2!} \QQ)  = \mc Hom(\pi_{1!} \QQ, \mc Hom(\pi_{2!} \QQ, \omega_S)) = \mc Hom(\pi_{1!} \QQ, \pi_{2*} D\QQ),$$
    using \cite[Proposition 2.6.3(ii)]{Kashiwara-Schapira} for the first equality. Identifying $D\QQ = \QQ[\dim X_2]$, since the sheaves $R^j\pi_{2*} \QQ$ are locally trivial we see from \cite[p.~624]{dC-Migliorini} that
    $$D(\pi_{1!} \QQ \otimes \pi_{2!} \QQ) =  \mc Hom(\pi_{1!} \QQ, \QQ[\dim S]) \otimes \pi_{2*} \QQ[\dim X_2 - \dim S] 
    $$
    $$= D\pi_{1!} \QQ \otimes \pi_{2 *} \QQ[\dim X_2 - \dim S] =  (\pi_{1*} \QQ \otimes \pi_{2*} \QQ)[\dim X_1 + \dim X_2 - \dim S].$$
    Thus the needed identity is obtained by applying Verdier duality to \eqref{eq:rel-Kuenneth} and shifting by $\dim X_1 \times_S X_2$. 

If $S$ is a smooth Deligne--Mumford stack, we choose an \'etale cover $p: S' \to S$ by a smooth variety $S'$. By smooth base change, the K\"unneth morphism $R\pi_{1*} \QQ \otimes R\pi_{2*} \QQ \to R\pi_* \QQ$ becomes an isomorphism after pulling back to $S'$. In particular, the mapping cone $C$ satisfies $p^* C = 0$, and thus $C = 0$.
\end{proof}

\begin{proof}[Proof of Corollary \ref{cor: splitting, GL}]
    By Proposition \ref{prop: relative B relations}, the character variety $M_B^{\GL}(\C_{g,1},n,d)$ is a quotient of $M_B^{\GL}(\C_{g,1},1,0) \times_{\M_{g,1}} M_B^{\SL}(\C_{g,1},n,d)$ by the relative group object $\widetilde\Gamma$ with fiber $\Gamma = \Hom(\pi_1(S_{g,1}), \mu_n)$ %
    over $\M_{g,1}$. %
	We have an isomorphism $R\pi_{\GL *} \QQ = (R\pi_{\GL_1 \times \SL} \QQ)^{\widetilde\Gamma}$, as in \eqref{eq:RpiGamma}.  
    By Lemma \ref{lemma: smooth relative Künneth}, using that $M_B^{\SL}(\C_{g,1},n,d)$ and $M_B^{\GL}(\C_{g,1},1,0)$ are locally trivial fibrations with smooth fibers, we see that
    $$R\pi_{\GL_1 \times \SL} \QQ = R\pi_{\GL_1 *} \QQ \otimes R\pi_{\SL *} \QQ.$$
    To show the splitting of $R\pi_{\GL_1 *} \QQ$, we use Simpson's ($d = 0$) nonabelian Hodge theorem
    $M_B^{\GL}(\C_{g,1},1,0) \cong M_\Dol^{\GL}(\C_{g,1},1,0)$
    \cite[Proposition~6.8 and Theorem~7.18]{Simpson-reps2}, together with the identification
    \[
    M_\Dol^{\GL}(\C_{g,1},1,0)
    \cong
    \Pic_{\C_{g,1}/\M_{g,1}}^0 \times_{\M_{g,1}} \mathbb E_{g,1},
    \]
    where $\mathbb E_{g,1}\to\M_{g,1}$ is the Hodge bundle with fiber $H^0(C,\omega_C)$ over $(C,p)$. This space retracts onto $\Pic_{\C_{g,1}/\M_{g,1}}^0$ over $\M_{g,1}$, so the splitting follows from Deligne's original theorem \cite[Proposition~2.1]{DeligneE2}.
    Furthermore, $\Gamma$ acts trivially on the cohomology of the fibers. This is because the action on the fibers can be extended to the multiplication action of the connected group $\Hom(\pi_1(S_{g,1}), \CC^*) \cong (\CC^*)^{2g}$ on itself. Then, since $M_B^{\GL}(\C_{g,1},1,0) \to \M_{g,1}$ is topologically locally trivial, we have $(R\pi_{\GL_1 *} \QQ)^{\widetilde\Gamma} = R\pi_{\GL_1 *} \QQ$. Thus, using \eqref{eq:RpiGamma}, we conclude that
    \begin{equation} \label{eq:RpiGL}
    \begin{aligned}
    R\pi_{\GL *} \QQ
    &\cong (R\pi_{\GL_1 \times \SL} \QQ)^{\widetilde\Gamma} \\
    &\cong (R\pi_{\GL_1 *} \QQ \otimes R\pi_{\SL *} \QQ)^{\widetilde\Gamma} \\
    &\cong R\pi_{\GL_1 *} \QQ \otimes (R\pi_{\SL *} \QQ)^{\widetilde\Gamma} \\
    &\cong R\pi_{\GL_1 *} \QQ \otimes R\pi_{\PGL *} \QQ.
    \end{aligned}
    \end{equation}
    By smooth base change \cite[Theorem~3.2.13(ii)]{Dimca}, Theorem~\ref{thm: splitting} gives a splitting
    \[
    R\pi_{\PGL *}\QQ \cong \bigoplus_i R^i\pi_{\PGL *}\QQ[-i]
    \]
    over $\M_{g,1}$. Combining this with the splitting of $R\pi_{\GL_1 *}\QQ$ above and \eqref{eq:RpiGL} gives the desired splitting of $R\pi_{\GL *}\QQ$.
\end{proof}

\section{Computations of the stable rational cohomology} \label{sec: computations}
In this section, we prove \Cref{thm:mainPGL,thm:sl-stability-answer} by combining the degeneration of the Leray--Serre spectral sequence at $E_2$ (Theorem \ref{thm: splitting} and Corollaries \ref{cor: splitting, GL}--\ref{cor: splitting, SL}) with previously known results on the cohomology of the fibers  and the cohomology of $\M_g$ and $\M_{g,1}$ with coefficients in symmetric and exterior powers of the symplectic local systems. In the smaller stable range $k \le s(g)$, feeding these local systems into Looijenga's computation of the cohomology of $\M_g$ with symplectic coefficients \cite{looijenga-symplectic-coeffs} yields Corollary~\ref{cor:stab-PGL}.

\subsection{Tautological generators and monodromy} \label{subsec: monodromy}

The universal $\PGL_n$-character variety $\M_B^{\PGL}(\C_g,n,d) \to \M_g$ is a topologically locally trivial fiber bundle, as are
\[
M_B^{\SL}(\C_{g,1},n,d) \to \M_{g,1} \quad \mbox{ and } \quad M_B^{\GL}(\C_{g,1},n,d) \to \M_{g,1}.
\]
Here, we identify the local systems $R^q \pi_* \QQ$ appearing in the Leray--Serre spectral sequences.  

As in \S\ref{sec: Betti}, we write $X_{g,n,d}$ for the subspace of $\Hom(\pi_1(S_g), \PGL_n(\CC))$ parametrizing maps whose obstruction to lifting to $\SL_n(\CC)$ is $d \in \ZZ/n\ZZ$. Fix a base point $p \in S_g$, and let $\widetilde S_g$ be the universal cover of $S_g$.

\begin{definition} \label{def: universal PGL bundle}
    The tautological principal $\PGL_n$-bundle over $S_g \times \M^{\PGL}_B(S_g,n,d)$ is $$ \mc P := \big(\widetilde S_g \times X_{g,n,d} \times \PGL_n(\CC)\big)/\big(\pi_1(S_g,p) \times \PGL_n(\CC)\big).$$
\end{definition} 
\noindent  
Here, the action of $\pi_1(S_g,p) \times \PGL_n(\CC)$ is given by
$$(\gamma, M) \cdot (x, \phi, N) := (\gamma x, M\phi M^{-1}, M \phi(\gamma) N).$$
The structure map to $S_g \times \M_B^{\PGL}(S_g,n,d)$ is induced by projection to the first two factors.

\begin{remark}
    The tautological principal $\PGL_n$-bundle on $S_g \times M_B^{\GL}(S_{g,1},n,d)$ is the pullback of $\mc P$ along the projectivization map $M_B^{\GL}(S_{g,1},n,d) \to \M_B^{\PGL}(S_g,n,d)$ \cite[\S 1]{hausel-thaddeus-generation-pgl2}, as is the tautological principal $\PGL_n$-bundle on $S_g \times M_B^{\SL}(S_{g,1},n,d)$ \cite[\S 4.1]{Hausel-RV}.
\end{remark}

We use the following convention for Chern classes of principal $\PGL_n$-bundles. The quotient map $\SL_n(\CC) \to \PGL_n(\CC)$ has finite kernel, and hence the induced pullback
\[
H^*(\BPGL_n(\CC);\QQ) \longrightarrow
H^*(\BSL_n(\CC);\QQ)
\]
is an isomorphism. For $2 \le j \le n$, let $c_j \in H^{2j}(\BPGL_n(\CC);\QQ)$ be the unique class whose pullback to $\BSL_n(\CC)$ is the usual Chern class of the universal rank $n$ bundle with trivial determinant. The generators for $H^*(\M_B^{\PGL}(S_g,n,d); \QQ)$ are the K\"unneth components of these classes. The formula below records this without choosing a
basis. By Poincar\'e duality on $S_g$, the operation
\[
a \longmapsto p_{2*}(a \cup p_1^*(x)), \qquad x \in H^i(S_g;\QQ),
\]
extracts the K\"unneth component of $a$ lying over $H^{2-i}(S_g;\QQ)$ and pairs it with $x$. 

Let $p_1$ and $p_2$ denote the projections of $S_g \times \M_B^{\PGL}(S_g,n,d)$ onto the first and second component, respectively.

\begin{definition}\label{def: alpha, beta, psi}
    For $i \in \{0,1,2\}$ and $j \in \ZZ_{\ge 2}$ we define the K\"unneth component $$V_{i,j} := \{p_{2*}(c_j(\mc P) \cup p_1^*(x)): x \in H^i(S_g;\QQ)\} \subset H^{i+2j-2}(\M_B^{\PGL}(S_g,n,d); \QQ).$$
\end{definition}

\begin{prop}\label{prop: free generation} %
    The subspaces $V_{i,j}$ for $0 \le i \le 2$ and $2 \le j \le n$ generate the ring $H^*(\M^{\PGL}_B(S_g,n,d); \QQ)$ with no relations in degrees less than or equal to $R_{g,n}$. Moreover, $\dim V_{i,j} = \dim H^i(S_g; \QQ)$.
\end{prop}

\begin{proof}
The generation statement is a consequence of Markman's tautological generation result for $M_{\Dol}^{\GL}(C,n,d)$ \cite[Theorem~7]{Markman}, as reinterpreted on the Betti side by Hausel and Thaddeus \cite[(5.1)]{hausel-thaddeus-generation-pgl2}. After choosing a homogeneous basis of $H^*(S_g; \QQ)$, there are no relations in degrees $ \le R_{g,n}$ among the images in $V_{*,j}$, by \cite[Lemma 4.1.12]{Hausel-RV}.
\end{proof}

Recall that $\mc P$ denotes the tautological principal $\PGL_n$-bundle on $S_g \times \M_B^{\PGL}(S_g,n,d)$.

\begin{prop} \label{prop: identifying monodromy}
The Chern classes $c_i(\mc P) \in H^{2i}(S_g \times \M_B^{\PGL}(S_g,n,d); \QQ)$ are fixed by the mapping class group action.
\end{prop}

\begin{proof}
    We follow the strategy of \cite[(5.1)]{hausel-thaddeus-generation-pgl2}. Represent the mapping class by a diffeomorphism $f$ of $S_g$ fixing the base point $p$; it acts on $X_{g,n,d}$ by $\phi \mapsto \phi \circ f_*^{-1}$. Both $\mc P$ and $f^* \mc P$ are $\PGL_n$-equivariant flat $\PGL_n$-bundles over $S_g \times X_{g,n,d}$, trivialized over $\{p\} \times X_{g,n,d}$, whose restriction to $S_g \times \{\phi\}$ realizes the monodromy representation $\phi$. A framed flat bundle is determined by its monodromy, so there is a unique $\PGL_n$-equivariant isomorphism $\mc P \xrightarrow{\sim} f^*\mc P$ preserving the framing. These isomorphisms are compatible with composition; in particular, $c_i(\mc P) = f^*c_i(\mc P)$.
\end{proof}

\begin{remark}\label{rmk: simplicial evaluation}
    Another proof of Proposition~\ref{prop: identifying monodromy}, following ideas from \cite{Shende}, goes as follows. Take a simplicial set with geometric realization homotopic to $S_g$, and view this as a simplicial scheme $X$ by assigning $\mathrm{Spec}\, \CC$ to each element. Here $\BPGL_n$ denotes the simplicial classifying scheme of $\PGL_n$. Then
    $$\M_B^{\PGL}(S_g,n,d) \cong \underline{\Hom}(X, \BPGL_n),$$
where $\underline{\Hom}$ is the internal Hom in the category of simplicial schemes. %
    From this perspective, the evaluation construction
    $$X \times {\underline{\Hom}}(X, \BPGL_n) \to \BPGL_n$$
    is functorial in the simplicial model $X$. Hence, self-homeomorphisms of $S_g$ preserve the pulled-back classes, which therefore define $\Mod_g$-invariant classes on the character stack. Passing to topological realization and forgetting the flat structure, this same evaluation map classifies the underlying topological principal bundle. Its adjoint is the comparison map $\varphi_{C,d}$ of \S\ref{subsec: comparison maps}. 
\end{remark}

We recall the notation established in the introduction \eqref{eq:Vg} and \eqref{eq: grading} and used in the statement of our main results. First,    $V_g \cong H^1(S_g;\QQ)$ viewed as local system on $\M_g$, i.e., the standard symplectic representation. Next,
\[
V_\circ = \QQ \oplus \QQ[-2], \quad \quad W_n = \QQ \oplus \QQ[-2] \oplus \cdots \oplus \QQ[2-2n], \quad \quad 
\quad  \widetilde{W}_n = W_n / \QQ. %
\]
Here, we think of $V_\circ$ as the even cohomology of $S_g$, viewed as a local system on $\M_g$ with trivial monodromy, supported in degrees 0 and 2. Likewise, $W_n$ and $\tW_n$ may be viewed as local systems on $\M_g$ with trivial monodromy, supported in even degrees. With this notation, we have the following.

\begin{prop} \label{prop: monodromy action}
Let $\pi \colon \M_B^{\PGL}(\C_g,n,d) \to \M_g$ be the projection. In degrees less than or equal to $R_{g,n}$, there is an isomorphism of graded local systems on $\M_g$
\[
R^\bullet \pi_* \QQ \;\cong\; \textstyle{\bigwedge^\bullet}(V_g \otimes \tW_n) \otimes \Sym^\bullet(V_\circ \otimes \tW_n).
\]
\end{prop}

\begin{proof}
Since $\M_B^{\PGL}(\C_g,n,d)\to\M_g$ is the bundle associated to the $\Mod_g$-action on the fiber $\M_B^{\PGL}(S_g,n,d)$ (Definition~\ref{defn: relative total PGL}), the local system $R^q\pi_*\QQ$ is the one associated to the graded $\Mod_g$-module $H^q(\M_B^{\PGL}(S_g,n,d);\QQ)$. It therefore suffices to identify this module in degrees at most $R_{g,n}$.
For $0\le i\le 2$ and $2\le j\le n$, Definition~\ref{def: alpha, beta, psi} gives a surjection
\[
H^i(S_g;\QQ)\longrightarrow V_{i,j},
\qquad
x\longmapsto p_{2*}(c_j(\mc P)\cup p_1^*x).
\]
It is $\Mod_g$-equivariant by Proposition~\ref{prop: identifying monodromy}, and is an isomorphism by Proposition~\ref{prop: free generation}. Hence, as graded $\Mod_g$-modules,
\[
V_{0,j}\cong \QQ[2-2j],\qquad
V_{1,j}\cong V_g[2-2j],\qquad
V_{2,j}\cong \QQ[-2j].
\]
Equivalently, as graded $\Mod_g$-modules, we have
\[
\bigoplus_{j=2}^n V_{1,j}\cong V_g\otimes \tW_n,
\qquad
\bigoplus_{j=2}^n(V_{0,j}\oplus V_{2,j})\cong V_\circ\otimes \tW_n.
\]
The first summand consists of odd-degree generators, and the second of even-degree generators. Proposition~\ref{prop: free generation} therefore gives
\[
H^*(\M_B^{\PGL}(S_g,n,d);\QQ)\cong
\textstyle{\bigwedge^\bullet}(V_g\otimes \tW_n)
\otimes
\Sym^\bullet(V_\circ\otimes \tW_n)
\]
in degrees less than or equal to $R_{g,n}$, as required. 
\end{proof}

\subsection{The Leray calculation} \label{subsec: leray}

We now proceed with the computation of the cohomology of universal character varieties in the range of degrees bounded above by $R_{g,n}$ and complete the proofs of \Cref{thm:mainPGL,thm:sl-stability-answer}.

\begin{proof}[Proof of Theorem \ref{thm:mainPGL}]
Consider the Leray--Serre spectral sequence
\[
E_2^{p,q} = H^p\big(\M_g;\, R^q\pi_*\QQ\big) \Longrightarrow H^{p+q}\big(\M_B^{\PGL}(\C_g,n,d);\QQ\big)
\]
of the fibration $\pi \colon \M_B^{\PGL}(\C_g,n,d) \to \M_g$. For $k \le R_{g,n}$, every term contributing to total degree $k$ has $q \le k \le R_{g,n}$, so by Proposition~\ref{prop: monodromy action} the degree-$k$ part of the $E_2$-page equals the degree-$k$ part of $H^*\big(\M_g;\, \bigwedge^\bullet(V_g \otimes \tW_n)\big) \otimes \Sym^\bullet(V_\circ \otimes \tW_n)$. By Theorem~\ref{thm: splitting} the spectral sequence degenerates at $E_2$, giving the theorem. 
\end{proof}

In our proof of Theorem~\ref{thm:sl-stability-answer}, we use Theorem~\ref{thm:mainPGL} together with the following proposition. Recall that the universal rank-one character variety $M_B^{\GL}(\C_{g,1},1,0)$ is pulled back from $\M_g$ (Remark~\ref{remark: point-pushing}). We write $M_B^{\GL}(S_g,1,0) := \Hom(\pi_1(S_g),\CC^*)$ for its fiber, on which $\Mod_g$ acts by precomposition.

\begin{prop}\label{prop: GL1 monodromy}
    The $\Mod_g$-representation $H^i(M_B^{\GL}(S_g,1,0); \QQ)$ is the $i^{\textrm{th}}$ exterior power of the standard symplectic representation.
\end{prop}

\begin{proof}
    We note that
    $$M_B^{\GL}(S_g,1,0) = \Hom(\pi_1(S_g), \CC^*) = \Hom(H_1(S_g; \ZZ), \CC^*) \cong (\CC^*)^{2g},$$
    so $H^i(M_B^{\GL}(S_g,1,0); \QQ) = \textstyle{\bigwedge^i} H^1(S_g;\QQ)$, and the proposition follows. 
\end{proof}

Let $\pi_G$ denote the projection map $M_B^G(\C_{g,1},n,d) \to \M_{g,1}$ for $G = \GL, \SL, \PGL$.     

\begin{proof}[Proof of Theorem \ref{thm:sl-stability-answer}]
Two identities of graded local systems on $\M_{g,1}$ drive the proof. In degrees $\le R_{g,n}$,
\[
R^\bullet\pi_{\SL *}\QQ = R^\bullet\pi_{\PGL *}\QQ,
\quad \mbox{ and } \quad 
R^\bullet\pi_{\GL *}\QQ = R^\bullet\pi_{\PGL *}\QQ \otimes \textstyle{\bigwedge^\bullet} V_g,
\]
where $R^\bullet\pi_{\PGL *}\QQ$ is understood as pulled back to $\M_{g,1}$. The first isomorphism follows from Lemma~\ref{lemma: SL invariance}. The second combines \eqref{eq:RpiGamma} and \eqref{eq:RpiGL} with $R^i\pi_{\GL_1 *}\QQ=\bigwedge^iV_g$ (Proposition~\ref{prop: GL1 monodromy}).

For $\SL$, the local systems $R^\bullet\pi_{\PGL *}\QQ$ are pulled back from those of Proposition~\ref{prop: monodromy action}, so the computation proving Theorem~\ref{thm:mainPGL} carries over to $\M_{g,1}$, with Corollary~\ref{cor: splitting, SL} supplying the degeneration, and gives the stated answer. For $\GL$, the second identity yields
\[
H^*(M_B^{\GL}(\C_{g,1},n,d); \QQ) \cong H^*\big(\M_{g,1};\, \textstyle{\bigwedge^\bullet}(V_g \otimes \tW_n) \otimes \Sym^\bullet(V_\circ \otimes \tW_n) \otimes \textstyle{\bigwedge^\bullet} V_g\big),
\]
matching the stated expression once we identify $(V_g \otimes \tW_n) \oplus V_g = V_g \otimes W_n$. 
\end{proof}

\subsection{Stable consequences} \label{subsec: stable consequences}

Next we deduce Corollaries~\ref{cor:stab-PGL} and \ref{cor:stab-G}, showing the stabilization of the cohomology of universal character varieties in the smaller range of degrees bounded above by $s(g) = \frac{2g-2}{3}$.

Let $S_{\langle \lambda \rangle}(V_g)$ be the irreducible representation of $\Sp(2g,\QQ)$ attached to the partition $\lambda$. When no confusion seems possible, we also write $S_{\langle \lambda \rangle}(V_g)$ for the induced local system on $\M_g$, viewed as a graded object in degree $0$.

\begin{lemma}\label{lemma: exterior decomposition}
In degrees less than or equal to $3g$, the graded local system $\bigwedge^\bullet(V_g \otimes \tW_n)$ on $\M_g$ is a direct sum of graded local systems $S_{\langle \lambda \rangle}(V_g)[-\ell]$. Every summand satisfies $|\lambda| \le \ell$ and $|\lambda| \equiv \ell \pmod 2$, and the summands are independent of $g$.
\end{lemma}

\begin{proof}
Since $\tW_n = \bigoplus_{i=1}^{n-1} \QQ[-2i]$, we have $V_g \otimes \tW_n = \bigoplus_{i=1}^{n-1} V_g[-2i]$, and therefore
\[
\textstyle{\bigwedge^\bullet}(V_g \otimes \tW_n) = \bigotimes_{i=1}^{n-1} \textstyle{\bigwedge^\bullet}\big(V_g[-2i]\big).
\]
Every generator lies in $V_g \otimes \tW_n$ and so has degree at least $3$. A monomial of degree less than or equal to $3g$ is thus a product of at most $g$ generators, and only the exterior powers $\textstyle{\bigwedge^r} V_g$ with $r \le g$ occur in such degrees. In this range the exterior power decomposes as:
\[
\textstyle{\bigwedge^r} V_g \cong
\begin{cases}
\bigoplus_{i=0}^{r/2} S_{\langle 1^{2i} \rangle}(V_g), & r \text{ even},\\
\bigoplus_{i=1}^{(r+1)/2} S_{\langle 1^{2i-1} \rangle}(V_g), & r \text{ odd},
\end{cases}
\qquad (r \le g),
\]
so each summand $S_{\langle 1^a \rangle}(V_g)$ of the $i$-th factor has $a \le r_i$ and $a \equiv r_i \pmod 2$. A summand of $\bigwedge^\bullet(V_g \otimes \tW_n)$ is a tensor product $\bigotimes_i \bigwedge^{r_i} V_g$ placed in degree $\ell = \sum_i r_i(2i+1)$; writing $r = \sum_i r_i$, we have $\ell \ge 3r$ and $\ell \equiv r \pmod 2$. The Littlewood--Richardson rule (see \cite{YoungDiagrammaticMethods}) decomposes the tensor product of these constituents into representations $S_{\langle \lambda \rangle}(V_g)$ with $|\lambda| \le r$ and $|\lambda| \equiv r \pmod 2$. Hence $|\lambda| \le r \le \ell$ and $|\lambda| \equiv \ell \pmod 2$. For $r \le g$ both the exterior decomposition and the Littlewood--Richardson coefficients are independent of $g$, and the lemma follows. \end{proof}

\begin{proof}[Proof of Corollaries~\ref{cor:stab-PGL} and \ref{cor:stab-G}]
We start with the proof of Corollary~\ref{cor:stab-PGL}. Fix $k\leq s(g)$. Theorem~\ref{thm:mainPGL} identifies $H^k(\M_B^{\PGL}(\C_g,n,d);\QQ)$ with the degree-$k$ part of
\[
H^*\big(\M_g;\textstyle{\bigwedge^\bullet}(V_g\otimes\tW_n)\big)\otimes\Sym^\bullet(V_\circ\otimes\tW_n).
\]
By Lemma~\ref{lemma: exterior decomposition}, every nonzero summand contributing to this degree has the form
\[
H^a\big(\M_g;S_{\langle\lambda\rangle}(V_g)\big),
\qquad
a=k-\ell-m,
\]
where $S_{\langle\lambda\rangle}(V_g)[-\ell]$ occurs in $\bigwedge^\bullet(V_g\otimes\tW_n)$, $|\lambda|\leq \ell$, $|\lambda|\equiv \ell\pmod 2$, and the integer $m$ is an even degree occurring in $\Sym^\bullet(V_\circ\otimes\tW_n)$. The  summands are independent of $g$. 

By \cite[Theorem~1.1]{looijenga-symplectic-coeffs}, together with the improved homological stability of \cite{boldsen,RW-resolutions}, the group $H^a(\M_g;S_{\langle\lambda\rangle}(V_g))$ has reached its stable value whenever $a+|\lambda|\leq s(g)$, and its stable value vanishes unless $a\equiv|\lambda|\pmod 2$. In the present situation,
\[
a+|\lambda|
=
k-\ell-m+|\lambda|
\leq k
\leq s(g),
\]
so every summand contributing to $H^k(\M_B^{\PGL}(\C_g,n,d);\QQ)$ is stable. Since its multiplicity is also independent of $g$, the same is true of the whole cohomology group. Moreover, if such a summand is nonzero, then $k=a+\ell+m$, which is congruent to $|\lambda| + |\lambda| \pmod 2$. Hence, the cohomology vanishes in odd degrees.

Finally, $\tW_{n+1}\cong\tW_n\oplus\QQ[-2n]$. The new summand first contributes in degree $2n+1$ to $\bigwedge^\bullet(V_g\otimes\tW_{n+1})$ and in degree $2n$ to $\Sym^\bullet(V_\circ\otimes\tW_{n+1})$. Hence the expression of Theorem~\ref{thm:mainPGL} is unchanged in degrees $k\leq2n-1$, which completes the proof of Corollary~\ref{cor:stab-PGL}.

The proof of Corollary~\ref{cor:stab-G} is the same, using Theorem~\ref{thm:sl-stability-answer} and the stable cohomology of $\M_{g,1}$ with symplectic coefficients \cite{Kawazumi,Randal-Williams}.
\end{proof}

\section{The stable cohomology ring} \label{sec: ring}

We prove Theorems~\ref{thm:stablePGL} and~\ref{thm:sl-stable-ring} by comparing universal character varieties with moduli spaces of surfaces in a background space, and using related scanning theorems from homotopy theory \cite{cohen2006surfaces,GMTW,GRW-monoids}. The argument for $\PGL_n$ was suggested to us by Oscar Randal-Williams and follows the method of \cite[Appendix~B]{Randal-Williams}. The Teichm\"uller models below may equally be viewed as moduli spaces of smooth complex curves with continuous maps to a background space; we use the topological term ``surfaces'' to match the language of the scanning results. Throughout this section, $BG$ denotes the topological classifying space of a complex Lie group $G(\CC)$. In particular, $\BPGL_n$ denotes $B(\PGL_n(\CC))$.

\subsection{Comparison with surfaces in $\BPGL_n$} \label{subsec: comparison maps}

For a simply connected target $Y$, obstruction theory and the orientation of $S_g$ identify $\pi_0\Map(S_g,Y)$ with $H_2(Y;\ZZ)$ via $f\mapsto f_*[S_g]$. For $\gamma\in H_2(Y;\ZZ)$, define
\[
\Sc_g(Y)_\gamma := \mc T_g \times_{\Mod_g} \Map(S_g,Y)_\gamma.
\]
This is the Teichm\"uller model of the component of the moduli space of maps $f\colon S_g\to Y$ with $f_*[S_g]=\gamma$. In the scanning theorem of Cohen--Madsen \cite{cohen2006surfaces}, the components are indexed by $\pi_2(Y)$. These two indexing conventions are equivalent because, for simply connected $Y$, the Hurewicz map $\pi_2(Y)\to H_2(Y;\ZZ)$ is an isomorphism.

Note that $\Sc_g(Y)_\gamma$ is the total space of a fibration over $\Sc_g(*)\simeq\M_g$ with fiber $\Map(S_g,Y)_\gamma$. It is weakly equivalent to the Borel construction model $E\Diff^+(S_g) \times_{\Diff^+(S_g)} \Map(S_g,Y)_\gamma$ for the oriented diffeomorphism group of $S_g$, whose identity component is contractible \cite{Earle-Eells}. Throughout this section, we pass freely between such Teichm\"uller models and Borel construction models of weakly equivalent moduli spaces of surfaces with a map to a background space.

The universal curve over $\M_B^{\PGL}(\C_g,n,d)$ has a tautological principal $\PGL_n$-bundle $\mc P$. After forgetting the flat structure, the fiberwise classifying maps of $\mc P$ determine a map
\[
\Phi_{g,d} \colon \M_B^{\PGL}(\C_g,n,d) \longrightarrow \Sc_g(\BPGL_n)_d
\]
over $\M_g$, whose restriction to the fiber over a point $[C]$ is
\[
\varphi_{C,d} \colon \M_B^{\PGL}(C,n,d) \longrightarrow \Map(C,\BPGL_n)_d,
\]
the adjoint of the classifying map of $\mc P$ on $C \times \M_B^{\PGL}(C,n,d)$.

Recall that $\M_B^{\PGL}(C,n,d)=[X_{g,n,d}/\PGL_n(\CC)]$. The framed flat bundle on $C \times X_{g,n,d}$ is $\PGL_n(\CC)$-equivariant for the conjugation action and therefore descends to $\mc P$ on the quotient stack; thus $\varphi_{C,d}$ is defined on the stack. Since $d$ is coprime to $n$, its effect on rational cohomology may be read on the coarse character variety through the isomorphism $H^*(\M_B^{\PGL}(C,n,d);\QQ) \cong H^*(M_B^{\PGL}(C,n,d);\QQ)$. The components of the target are indexed by $H_2(\BPGL_n;\ZZ)=\ZZ/n\ZZ$, and the image lies in the component indexed by $d$.

\begin{prop} \label{prop: fiber comparison}
For a fixed curve $C$ of genus $g$, the map $\varphi_{C,d}$ induces an isomorphism on rational cohomology in degrees $\le R_{g,n}$.
\end{prop}

\begin{proof}
Let $\ev \colon C \times \Map(C,\BPGL_n)_d \to \BPGL_n$ be the evaluation map. Since $\varphi_{C,d}$ is adjoint to the classifying map of $\mc P$, we have
\[
(\id_C \times \varphi_{C,d})^*\ev^*c_j=c_j(\mc P).
\]
Thus, $\varphi_{C,d}^*$ sends the K\"unneth components of $\ev^*c_j$ to the subspaces $V_{i,j}$ of Definition~\ref{def: alpha, beta, psi}.

It remains to compare generators and relations. The classes $c_2,\dots,c_n$ give a rational equivalence $\BPGL_n \xrightarrow{\sim_\QQ} \prod_{j=2}^n K(\QQ,2j)$. Haefliger's mapping-space calculation \cite{Haefliger} shows that the induced map of mapping spaces is a rational equivalence on every component; in particular, the finite component label $d\in\ZZ/n\ZZ$ disappears after rationalization. For each factor $K(\QQ,2j)$, \cite[Lemma~3.5.3]{BDPW} computes the pointed mapping space, while the (split) evaluation fibration supplies the additional factor $K(\QQ,2j)$. Thus, for every $d$,
\[
\Map(C,\BPGL_n)_d \simeq_{\QQ} \prod_{j=2}^n \Big(K(\QQ,2j) \times K\big(H^1(C;\QQ),2j-1\big) \times K(\QQ,2j-2)\Big).
\]
The right side is independent of $d$, so $H^*(\Map(C,\BPGL_n)_d;\QQ)$ has the same free generators for every component. These generators map to $V_{i,j}$, whose description in \Cref{prop: free generation} is likewise uniform for $d\in(\ZZ/n\ZZ)^*$. The proposition follows.
\end{proof}

\begin{lemma} \label{lem: range comparison}
Let $f \colon E \to E'$ be a map over a space $B$, and write $p \colon E \to B$ and $p' \colon E' \to B$ for the structure maps. If the induced morphism $R^q p'_*\QQ \to R^q p_*\QQ$ 
is an isomorphism for $q \le N$, then $f^* \colon H^k(E';\QQ) \to H^k(E;\QQ)$ 
is an isomorphism for $k \le N$.
\end{lemma}

\begin{proof}
The hypothesis says that the induced map of canonical truncations
\[
\tau_{\le N}Rp'_*\QQ \longrightarrow \tau_{\le N}Rp_*\QQ
\]
is an isomorphism. For either derived direct image, the cone of the canonical map from its truncation to the full complex lies in $D^{\geq N+1}(B)$, so its hypercohomology vanishes in degrees $\le N$. Thus, in this range, the cohomology of $E'$ and $E$ is computed by $R\Gamma(B,\tau_{\le N}Rp'_*\QQ)$ and $R\Gamma(B,\tau_{\le N}Rp_*\QQ)$, respectively. The asserted isomorphism follows.
\end{proof}

\begin{prop} \label{prop: family comparison}
The map $\Phi_{g,d}$ induces an isomorphism on rational cohomology in degrees $\le R_{g,n}$.
\end{prop}

\begin{proof}
Write $p_{\Sc} \colon \Sc_g(\BPGL_n)_d \to \M_g$ and $\pi \colon \M_B^{\PGL}(\C_g,n,d) \to \M_g$ for the projections. The restriction of $\Phi_{g,d}$ over $[C]$ is $\varphi_{C,d}$, so \Cref{prop: fiber comparison} gives an isomorphism 
\[
R^q(p_{\Sc})_*\QQ \xrightarrow{\ \sim\ } R^q\pi_*\QQ, \qquad \mbox{ for \ } q \le R_{g,n}.
\]
The result then follows from \Cref{lem: range comparison}.
\end{proof}

\subsection{Scanning and the $\PGL_n$ ring} \label{subsec: scanning PGL ring}

We use the following finite-genus form of the scanning theorem for surfaces in a background space.

\begin{theorem} \label{thm: scanning}
Let $Y$ be a simply connected CW complex of finite rational type, and let $\gamma \in H_2(Y;\ZZ)$. After translating the component met by the natural scanning map to the zero component, there is a scanning map
\[
\alpha_{g,\gamma} \colon \Sc_g(Y)_\gamma \longrightarrow \Omega^\infty_0\big(MTSO(2) \wedge Y_+\big)
\]
which induces an isomorphism on rational cohomology in degrees $\le s(g)$.
\end{theorem}

\noindent Here, $Y_+$ denotes $Y$ with a disjoint basepoint, $\wedge$ is the smash product, and $\Omega^\infty_0E$ is the zero (basepoint) component of the infinite loop space $\Omega^\infty E$. Translation to this component is possible because the infinite loop space is grouplike. Also, $MTSO(2)$ is the Thom spectrum of the inverse of the tautological complex line bundle $L \to \BSO(2)=\CC\PP^\infty$. This finite-genus scanning theorem is due to Cohen--Madsen \cite[Theorems~0.1 and~0.3]{cohen2006surfaces}, whose stable identification rests on the theorem of Galatius--Madsen--Tillmann--Weiss on the classifying space of the cobordism category of surfaces mapping to $Y$ \cite[\S 5]{GMTW}. In the stated range, it follows by combining their stable identification, or its formulation for tangential structures \cite[Corollary~F]{GRW-monoids}, with homological stability for surfaces in a background space \cite[\S 7.4]{RW-resolutions}, cf. \cite[Theorem~4.17]{boldsen}.

For a spectrum $E$ of finite rational type, one has $\pi_k(\Omega^\infty_0E)\otimes\QQ \cong H_k(E;\QQ)$, for $k > 0$. Hence $\Omega^\infty_0E$ is rationally a product of Eilenberg--MacLane spaces, and
\begin{equation} \label{eq: rational infinite loop cohomology}
H^*(\Omega^\infty_0E;\QQ) \cong \Sym^\bullet\big(H^{>0}(E;\QQ)\big)
\end{equation}
as a free graded-commutative algebra. In \eqref{eq: rational infinite loop cohomology}, the inclusion of the generating space $H^{>0}(E;\QQ)$ is the stable cohomology suspension
\[
H^{>0}(E;\QQ)\longrightarrow H^{>0}(\Omega^\infty_0E;\QQ).
\]
Applied to $E=MTSO(2)\wedge Y_+$, the Thom isomorphism gives
\begin{equation} \label{eq: scanning algebra general}
H^*\big(\Omega^\infty_0(MTSO(2)\wedge Y_+);\QQ\big)
\cong
\Sym^\bullet\!\left(\big(\QQ[\psi]\otimes H^*(Y;\QQ)[2]\big)_{>0}\right),
\end{equation}
where $e:=c_1(L)$ and $\psi:=-e$ is the first Chern class of the relative cotangent line bundle.

We now take $Y=\BPGL_n$ and restrict to the component $\Sc_g(\BPGL_n)_d$. The positive-degree classes obtained by stable cohomology suspension in \eqref{eq: rational infinite loop cohomology} are primitive and hence fixed by translation. Thus the resulting algebra and its distinguished generators are independent of $d$. Since $H^*(\BPGL_n;\QQ)=\QQ[c_2,\dots,c_n]$, \eqref{eq: scanning algebra general} becomes
\begin{equation} \label{eq: scanning algebra PGL}
H^*\big(\Omega^\infty_0(MTSO(2)\wedge(\BPGL_n)_+);\QQ\big)
\cong \Sym^\bullet(K_n).
\end{equation}
Under scanning, the generator indexed by a monomial $\psi^{a_0}c_2^{a_2}\cdots c_n^{a_n}$ becomes the generalized $\kappa$-class $\varpi_!\big(\psi^{a_0}u^*c_2^{a_2}\cdots u^*c_n^{a_n}\big)$ 
on $\Sc_g(\BPGL_n)_d$, where $\varpi\colon\mathcal E\to\Sc_g(\BPGL_n)_d$ is the universal oriented surface bundle and $u\colon\mathcal E\to\BPGL_n$ is the universal map.

\begin{proof}[Proof of \Cref{thm:stablePGL}]
The preceding description of scanning and the definition of $\Phi_{g,d}$ identify the tautological homomorphism \eqref{eq: PGLring} with the composite
\[
\Sym^\bullet(K_n)
\xrightarrow{\ \alpha_{g,d}^*\ }
H^*\big(\Sc_g(\BPGL_n)_d;\QQ\big)
\xrightarrow{\ \Phi_{g,d}^*\ }
H^*\big(\M_B^{\PGL}(\C_g,n,d);\QQ\big).
\]
By \Cref{thm: scanning} and \Cref{prop: family comparison}, the two maps are isomorphisms in degrees $\le s(g)$ and $\le R_{g,n}$, respectively, and the theorem follows. The source algebra, its generators, and both comparisons are uniform in $d\in(\ZZ/n\ZZ)^*$.
\end{proof}

\subsection{Surfaces with one marked point} \label{subsec: marked point}

Let $Y$ be a simply connected CW complex of finite rational type and let $\gamma \in H_2(Y;\ZZ)$. Write $\Sc_{g,1}(Y)_\gamma := \mc T_{g,1} \times_{\Mod_{g,1}} \Map(S_g,Y)_\gamma$ for the moduli space of genus-$g$ surfaces with a marked point equipped with a map to $Y$ in the component $\gamma$. The forgetful map $\M_{g,1}\to\M_g$ is the universal curve. Equivalently, by the Birman exact sequence \cite[\S 4.2]{Farb-Margalit}, the homotopy fiber of $B\Mod_{g,1}\to B\Mod_g$ is $B\pi_1(S_g)\simeq S_g$. Applying the same Borel construction with $\Map(S_g,Y)_\gamma$ identifies $\Sc_{g,1}(Y)_\gamma$ with the total space of the universal surface bundle $\varpi \colon \mathcal E \to \Sc_g(Y)_\gamma$, equipped with its universal map $u \colon \mathcal E \to Y$; under this identification, $u$ is evaluation at the marked point. We retain the notation $\psi := c_1(T_\varpi^\vee)$.

\begin{prop} \label{prop: marked point}
Let $Y$ be a simply connected CW complex of finite rational type, and let $\gamma \in H_2(Y;\ZZ)$. The algebra map
\[
\Theta \colon H^*\big(\Sc_g(Y)_\gamma;\QQ\big) \otimes \QQ[\psi] \otimes H^*(Y;\QQ) \longrightarrow H^*\big(\Sc_{g,1}(Y)_\gamma;\QQ\big), \qquad \alpha \otimes \psi^a \otimes c \longmapsto \varpi^*\alpha \cdot \psi^a \cdot u^*c,
\]
is an isomorphism in degrees less than or equal to $s(g)$.
\end{prop}

The proof uses the boundary form of scanning, a Leray--Hirsch argument through degree $s(g)$, and the Gysin sequence of the circle bundle of tangent rays. The scanning results we use are formulated for surfaces with a parametrized boundary component, whereas \Cref{prop: marked point} concerns surfaces with a marked point. To pass between these settings, we equip the marked point with a tangent ray and remove a small disc around it. The next lemma makes this comparison precise.

\begin{lemma} \label{lem: framed point boundary}
Let $Y$ be a simply connected CW complex of finite rational type and let $\gamma \in  H_2(Y;\ZZ)$. Let $\Sc^{\mathrm{fr}}_g(Y)_\gamma$ denote the moduli space of genus-$g$ surfaces with a marked point, a tangent ray at that point, and a map to $Y$ in the component $\gamma$, and let
\[
\varrho \colon \Sc^{\mathrm{fr}}_g(Y)_\gamma \longrightarrow \Sc_{g,1}(Y)_\gamma
\]
forget the tangent ray. This is the circle bundle of tangent rays in the vertical tangent bundle. The evaluation map $\ev \colon \Sc^{\mathrm{fr}}_g(Y)_\gamma \longrightarrow Y$ is a fibration. For $y_0\in Y$, its fiber is weakly equivalent to the moduli space $\Sc^{\partial}_g(Y)_\gamma$ of genus-$g$ surfaces with one parametrized boundary component and a map to $Y$ that is constant equal to $y_0$ on the boundary.
\end{lemma}

\begin{proof}
Write $\Diff_g^+:=\Diff^+(S_g)$, and let $\Diff^+_{g;p,v}\subset \Diff_g^+$ be the subgroup fixing a point $p$ and a tangent ray $v$ at $p$. Since $\Diff_g^+$ acts transitively on tangent rays with stabilizer $\Diff^+_{g;p,v}$, there is a weak equivalence
\[
\Sc^{\mathrm{fr}}_g(Y)_\gamma\simeq E\Diff_g^+\times_{\Diff^+_{g;p,v}}\Map(S_g,Y)_\gamma.
\]
Evaluation at $p$ is a $\Diff^+_{g;p,v}$-equivariant fibration $\Map(S_g,Y)_\gamma\to Y$, where $\Diff^+_{g;p,v}$ acts trivially on $Y$, and therefore descends to the asserted fibration. Its fiber over $y_0$ is
\[
E\Diff_g^+\times_{\Diff^+_{g;p,v}}\Map_{y_0},
\]
where $\Map_{y_0}$ consists of maps taking $p$ to $y_0$.

Fix a closed disc $D\subset S_g$ centered at $p$ whose chosen parametrization carries the positive horizontal ray to $v$. Let $\Diff^+_{g;D}\subset \Diff^+_{g;p,v}$ be the subgroup fixing $D$ pointwise, and let $\Map_D\subset\Map_{y_0}$ consist of maps that are constant equal to $y_0$ on $D$. Restriction to $D$ and evaluation at $p$ form a map of Serre fibrations whose map on bases $\Map(D,Y)\to Y$ is a homotopy equivalence. Hence $\Map_D\to\Map_{y_0}$ is a $\Diff^+_{g;D}$-equivariant weak equivalence. By the isotopy extension theorem, restriction to $D$ gives a fibration
\[
\Diff^+_{g;D}\longrightarrow \Diff^+_{g;p,v}\longrightarrow \operatorname{Emb}^+_{p,v}(D,S_g),
\]
where $\operatorname{Emb}^+_{p,v}(D,S_g)$ is the space of orientation-preserving parametrized discs centered at $p$ whose first tangent vector spans the ray $v$. The derivative-at-the-center map $\operatorname{Emb}^+(D,S_g)\to\operatorname{Fr}^+(TS_g)$ is a fibration with contractible fibers \cite[proof of Proposition~11.2]{RW-resolutions}. Its inverse image over the contractible subspace of oriented frames based at $p$ whose first vector spans $v$ is therefore contractible. Thus $\operatorname{Emb}^+_{p,v}(D,S_g)$ is contractible and $\Diff^+_{g;D}\to \Diff^+_{g;p,v}$ is a weak equivalence. Passing to homotopy quotients gives
\[
\ev^{-1}(y_0)\simeq E\Diff_g^+\times_{\Diff^+_{g;D}}\Map_D=: \Sc^{\partial}_g(Y)_\gamma.
\]
Since $Y$ is simply connected, forgetting the condition that $p$ map to $y_0$ induces a bijection from pointed to unpointed homotopy classes of maps $S_g\to Y$, so this fiber lies in a single component indexed by $\gamma$. Removing the interior of $D$ identifies it with the Cohen--Madsen moduli space $S_{g,1}(Y;\mathrm{const})$ of surfaces with one parametrized boundary component mapped constantly to $y_0$.
\end{proof}

\begin{lemma} \label{lem: Leray-Hirsch}
Let $F \to E \to B$ be a fibration with $B$ simply connected and of finite rational type. Suppose that a finite set of homogeneous classes of positive degree in $H^*(E;\QQ)$ restrict to free generators of $H^*(F;\QQ)$, as a graded-commutative $\QQ$-algebra, through degree $N$. Let $A$ be the free graded-commutative $\QQ$-algebra on corresponding formal variables. Then cup product induces an isomorphism through degree $N$
\[
H^*(B;\QQ)\otimes A\longrightarrow H^*(E;\QQ).
\]
\end{lemma}

\begin{proof}
The Serre spectral sequence of the fibration has $E_2^{p,q} = H^p(B;\QQ) \otimes H^q(F;\QQ)$. For $q \le N$ the monomials in the formal variables restrict to a basis of $H^q(F;\QQ)$, so the restriction $H^q(E;\QQ) \to H^q(F;\QQ)$ is surjective and $E_\infty^{0,q} = E_2^{0,q}$: every differential vanishes on the fiber column in these degrees. Since $E_2^{p,q} = E_2^{p,0} \cdot E_2^{0,q}$ and the differentials are derivations vanishing on the base row, they vanish wherever the fiber degree is at most $N$, on every page. A differential whose source or target has total degree $\le N$ has source of fiber degree at most $N$, hence vanishes; so $E_\infty^{p,q} = E_2^{p,q}$ for $p+q \le N$. The cup-product map respects the filtrations, and on associated graded pieces in total degrees $\le N$ it induces the isomorphisms $H^p(B;\QQ) \otimes A_q \xrightarrow{\ \sim\ } H^p(B;\QQ) \otimes H^q(F;\QQ)$, where $A_q$ is the degree $q$ part of $A$.
\end{proof}

\begin{proof}[Proof of \Cref{prop: marked point}]
For $c\in H^*(Y;\QQ)$ homogeneous and $a\geq0$, write
\[
\kappa_a(c):=\varpi_!\big(\psi^a u^*c\big).
\]
By \Cref{thm: scanning} and \eqref{eq: scanning algebra general}, the classes $\kappa_a(c)$ freely generate $H^*(\Sc_g(Y)_\gamma;\QQ)$ as a graded-commutative $\QQ$-algebra through degree $s(g)$. Here $c$ ranges over a homogeneous basis of $H^*(Y;\QQ)$, $|c|$ denotes the degree of $c$, and  $2a+|c|>2$. Since $Y$ is of finite rational type, only finitely many generators have degree $\le s(g)$.

\smallskip
\noindent\emph{Step 1: framed points and boundary.} By \Cref{lem: framed point boundary}, the evaluation fibration
\[
\ev=u\circ\varrho\colon\Sc^{\mathrm{fr}}_g(Y)_\gamma\longrightarrow Y
\]
has fiber weakly equivalent to $\Sc^{\partial}_g(Y)_\gamma$. The boundary form of the scanning theorem \cite[Theorem~0.1]{cohen2006surfaces}, together with the homological stability used in \Cref{thm: scanning}, identifies the rational cohomology of this fiber through degree $s(g)$ with the free graded-commutative algebra on the restrictions of the classes $\kappa_a(c)$.

For a surface bundle with parametrized boundary, the boundary generalized $\kappa$-classes are defined by the same fiber-integral formula; these are the classes appearing in the stable calculation of \cite[Corollary~0.2]{cohen2006surfaces}. The restrictions above agree with these boundary classes. Indeed, capping the boundary by the parametrized disc contributes nothing to the fiber integral: the map to $Y$ is constant on the disc, so $u^*c$ vanishes there when $|c|>0$, while the parametrization trivializes the vertical tangent bundle, so $\psi$ vanishes there when $a>0$. The only case not covered by these vanishing statements is $a=|c|=0$, which is excluded by $2a+|c|>2$.

\smallskip

\noindent\emph{Step 2: evaluation and Leray--Hirsch.} Pulling back the finitely many classes $\kappa_a(c)$ of degree $\le s(g)$ along $\varpi\circ\varrho$ gives classes on $\Sc^{\mathrm{fr}}_g(Y)_\gamma$ whose restrictions freely generate the fiber cohomology through degree $s(g)$. Applying \Cref{lem: Leray-Hirsch} and using the preceding scanning identification gives, in degrees $\le s(g)$, an isomorphism
\begin{equation} \label{eq: marked framed Leray-Hirsch}
\overline\Theta\colon H^*\big(\Sc_g(Y)_\gamma;\QQ\big)\otimes H^*(Y;\QQ)\longrightarrow H^*\big(\Sc^{\mathrm{fr}}_g(Y)_\gamma;\QQ\big),\qquad \alpha\otimes c\longmapsto\varrho^*\varpi^*\alpha\cdot\ev^*c.
\end{equation}
Since $\ev=u\circ\varrho$, the map $\overline\Theta$ factors through $\varrho^*$: explicitly, $\overline\Theta(\alpha\otimes c)=\varrho^*\big(\varpi^*\alpha\cdot u^*c\big)$ in this range. Since $\overline\Theta$ is an isomorphism, it follows that $\varrho^*$ is surjective through degree $s(g)$.

\smallskip

\noindent\emph{Step 3: the Gysin sequence.} The Euler class of the circle bundle $\varrho$ is $e(T_\varpi)=-\psi$. The Gysin sequence and the surjectivity of $\varrho^*$ give short exact sequences
\[
0\longrightarrow H^{k-2}\big(\Sc_{g,1}(Y)_\gamma;\QQ\big)\xrightarrow{\ \cup\,\psi\ }H^k\big(\Sc_{g,1}(Y)_\gamma;\QQ\big)\xrightarrow{\ \varrho^*\ }H^k\big(\Sc^{\mathrm{fr}}_g(Y)_\gamma;\QQ\big)\longrightarrow0,
\qquad k\leq s(g).
\]
Because $\varrho^*\psi=0$, the composite $\varrho^*\circ\Theta$ factors through the quotient of the source of $\Theta$ by the ideal $(\psi)$. Under the natural identification of this quotient with $H^*(\Sc_g(Y)_\gamma;\QQ)\otimes H^*(Y;\QQ)$, the induced map is the isomorphism \eqref{eq: marked framed Leray-Hirsch}.

We prove that $\Theta$ is an isomorphism by induction on degree, with negative degrees understood to vanish. Fix $k\leq s(g)$ and suppose the result holds in degrees less than $k$. For surjectivity, let $z\in H^k(\Sc_{g,1}(Y)_\gamma;\QQ)$ and choose $w_0\in H^*(\Sc_g(Y)_\gamma;\QQ)\otimes H^*(Y;\QQ)$ such that $\overline\Theta(w_0)=\varrho^*z$. Regard $w_0$ as an element of the source of $\Theta$ of $\psi$-degree zero. Exactness gives
\[
z-\Theta(w_0)=\psi y
\]
for some $y\in H^{k-2}(\Sc_{g,1}(Y)_\gamma;\QQ)$. By induction, $y=\Theta(w_1)$ for some $w_1$, and hence $z=\Theta(w_0+\psi w_1)$. For injectivity, let $w$ be a degree-$k$ element in the kernel of $\Theta$. Then $w=\psi w_1$ for some $w_1$ of degree $k-2$. Since cup product with $\psi$ is injective, it follows that $\Theta(w_1)=0$. By induction, $w_1=0$, and hence $w=0$.
\end{proof}

\subsection{The special linear case} \label{subsec: SL ring}

Let $q\colon M_B^{\SL}(\C_{g,1},n,d)\to
\M_B^{\PGL}(\C_{g,1},n,d)$ 
be the quotient by $\widetilde\Gamma$ of \Cref{prop: relative B relations}. Base change of $\Phi_{g,d}$ along $\M_{g,1}\to\M_g$ gives a Cartesian square
\[
\begin{tikzcd}
\M_B^{\PGL}(\C_{g,1},n,d)
    \arrow[r,"\Phi^{(1)}_{g,d}"]
    \arrow[d]
&
\Sc_{g,1}(\BPGL_n)_d
    \arrow[d]
\\
\M_B^{\PGL}(\C_g,n,d)
    \arrow[r,"\Phi_{g,d}"]
&
\Sc_g(\BPGL_n)_d .
\end{tikzcd}
\]
The restriction of $\Phi^{(1)}_{g,d}$ to the fiber over $(C,p)$ is the map $\varphi_{C,d}$ of \Cref{prop: fiber comparison}.

Both $q$ and $\Phi^{(1)}_{g,d}$ are maps over $\M_{g,1}$. By Lemmas~\ref{lemma: SL invariance} and~\ref{lem: range comparison}, pullback along $q$ is an isomorphism on rational cohomology in degrees at most $R_{g,n}$. Likewise, Proposition~\ref{prop: fiber comparison} and Lemma~\ref{lem: range comparison} show that pullback along $\Phi^{(1)}_{g,d}$ is an isomorphism in the same range.

\begin{proof}[Proof of \Cref{thm:sl-stable-ring} for $G=\SL_n$]
Consider the composite
\[
\begin{aligned}
\Sym^\bullet(K_n)\otimes\QQ[\psi,c_2,\dots,c_n]
&\xrightarrow{\ \alpha_{g,d}^*\otimes 1\ }
H^*\big(\Sc_g(\BPGL_n)_d;\QQ\big)
    \otimes\QQ[\psi]\otimes H^*(\BPGL_n;\QQ)
\\
&\xrightarrow{\ \Theta\ }
H^*\big(\Sc_{g,1}(\BPGL_n)_d;\QQ\big)
\\
&\xrightarrow{\ \Phi_{g,d}^{(1)*}\ }
H^*\big(\M_B^{\PGL}(\C_{g,1},n,d);\QQ\big)
\\
&\xrightarrow{\ q^*\ }
H^*\big(M_B^{\SL}(\C_{g,1},n,d);\QQ\big).
\end{aligned}
\]
Here we use $H^*(\BPGL_n;\QQ)=\QQ[c_2,\dots,c_n]$. The first arrow is an isomorphism through degree $s(g)$ by \Cref{thm: scanning} and \eqref{eq: scanning algebra PGL}, and the second is an isomorphism through degree $s(g)$ by \Cref{prop: marked point}. The last two arrows are isomorphisms through degree $R_{g,n}$ by the preceding paragraph. Hence, the composite is an isomorphism in degrees at most $s(g)$.

It remains to identify the composite with the map \eqref{eq: SLring}. For $m=\psi^{a_0}c_2^{a_2}\cdots c_n^{a_n}$,
write
\[
\kappa(m):=
\varpi_!\left(
\psi^{a_0}u^*c_2^{a_2}\cdots u^*c_n^{a_n}
\right)
\in H^*\big(\Sc_g(\BPGL_n)_d;\QQ\big).
\]
The Cartesian squares of universal curves, the compatibility of their marked sections and classifying maps, and the base-change formula for Gysin maps give
\[
q^*\Phi_{g,d}^{(1)*}\varpi^*\kappa(m)
=
(p_{\SL_n})_!\left(
\psi_{\SL_n}^{a_0}
c_2(\mathcal P_{\SL_n})^{a_2}\cdots
c_n(\mathcal P_{\SL_n})^{a_n}
\right),
\]
together with
\[
q^*\Phi_{g,d}^{(1)*}\psi
=
\sigma_{\SL_n}^*\psi_{\SL_n}
\quad \mbox{ and }
\quad 
q^*\Phi_{g,d}^{(1)*}u^*c_j
=
\sigma_{\SL_n}^*c_j(\mathcal P_{\SL_n}),
\ \  \mbox{ for \ } 2\leq j\leq n.
\]
These are the images prescribed in \eqref{eq: SLring}, and the theorem follows.
\end{proof}

\subsection{The general linear case} \label{subsec: GL ring}

A local system parametrized by $M_B^{\GL}(\C_{g,1},n,d)$ has monodromy $e^{2\pi i d/n}\Id_n$ around the puncture and therefore does not extend to a local system on the closed curve. However, its projectivization and determinant do extend: the scalar $e^{2\pi i d/n}\Id_n$ maps to the identity in $\PGL_n(\CC)$, while its determinant is $e^{2\pi i d}=1$. Moreover, 
\[
(\det,\mathrm{proj})\colon \GL_n(\CC)\longrightarrow \CC^*\times\PGL_n(\CC)
\]
is surjective with finite kernel $\mu_n$, so the induced map
\[
\BGL_n\longrightarrow B\CC^*\times\BPGL_n
\]
is a rational cohomology equivalence. We therefore compare the general linear character variety with surfaces in the product background space $B\CC^*\times\BPGL_n$.

Let $J:=M_B^{\GL}(\C_{g,1},1,0)$, and set
\[
\Sc^\sigma_{g,1}(B\CC^*)_0
:=
\mc T_{g,1}\times_{\Mod_{g,1}}\Map_*(S_g,B\CC^*)_0.
\]
Thus $\Sc^\sigma_{g,1}(B\CC^*)_0$ parametrizes a marked surface together with a degree-zero map to $B\CC^*$ that sends the marked point to the basepoint. The tautological flat line bundle on the universal curve over $J$, with its trivialization along the marked section, has a based classifying map and hence induces a map
\[
\lambda\colon J\longrightarrow\Sc^\sigma_{g,1}(B\CC^*)_0.
\]

Set
\[
T_{g,n,d}:=
\Sc^\sigma_{g,1}(B\CC^*)_0
\times_{\M_{g,1}}
\Sc_{g,1}(\BPGL_n)_d.
\]
Let $\delta\colon M_B^{\GL}(\C_{g,1},n,d)\to J$ be induced by the determinant, and let $\mathrm{proj}\colon M_B^{\GL}(\C_{g,1},n,d)\to\M_B^{\PGL}(\C_{g,1},n,d)$ be induced by projectivization. The determinant has trivial monodromy around the puncture, so it defines a degree-zero local system on the closed curve and hence the map $\delta$. The pair of maps
\[
\begin{aligned}
\lambda\circ\delta&\colon M_B^{\GL}(\C_{g,1},n,d)\longrightarrow\Sc^\sigma_{g,1}(B\CC^*)_0,\\
\Phi^{(1)}_{g,d}\circ\mathrm{proj}&\colon M_B^{\GL}(\C_{g,1},n,d)\longrightarrow\Sc_{g,1}(\BPGL_n)_d
\end{aligned}
\]
together define a map over $\M_{g,1}$:
\[
\Psi_{\GL}\colon M_B^{\GL}(\C_{g,1},n,d)\longrightarrow T_{g,n,d}.
\]

\begin{proof}[Proof of \Cref{thm:sl-stable-ring} for $G=\GL_n$]
We first show that $\Psi_{\GL}^*$ is an isomorphism on rational cohomology through degree $R_{g,n}$. Write $\pi_J$, $\pi_{\SL}$, and $\pi_{\PGL}$ for the projections to $\M_{g,1}$. By the quotient description in Proposition~\ref{prop: relative B relations}, the relative K\"unneth isomorphism of Lemma~\ref{lemma: smooth relative Künneth}, and exactness of invariants with rational coefficients, we have
\begin{align*}
R\pi_{\GL *}\QQ
&\cong \big(R\pi_{J *}\QQ\otimes R\pi_{\SL *}\QQ\big)^{\widetilde\Gamma}\\
&\cong R\pi_{J *}\QQ\otimes \big(R\pi_{\SL *}\QQ\big)^{\widetilde\Gamma}\\
&\cong R\pi_{J *}\QQ\otimes R\pi_{\PGL *}\QQ,
\end{align*}
as in \eqref{eq:RpiGL}. This uses only the identification of the invariant part of the $\SL_n$ cohomology with the $\PGL_n$ cohomology. The resulting isomorphism is induced by determinant and projectivization, up to the pullback by multiplication by $n$ on the $J$-factor: indeed, after pullback to $J\times_{\M_{g,1}}M_B^{\SL}(\C_{g,1},n,d)$, the determinant is the $n$-th power of the tautological rank-one local system. Multiplication by $n$ is a rational cohomology equivalence on every fiber of $J\to\M_{g,1}$, so determinant and projectivization induce an isomorphism on the cohomology local systems in every degree.

On each fiber, the restriction of $\lambda$ is the map
\[
\Hom(\pi_1(S_g),\CC^*)\longrightarrow\Map_*(S_g,B\CC^*)_0,
\]
which is a weak equivalence because both spaces are $K(H^1(S_g;\ZZ),1)$ and the tautological classifying map induces the standard identification on fundamental groups. The restriction of $\Phi^{(1)}_{g,d}$ induces an isomorphism on rational cohomology through degree $R_{g,n}$ by \Cref{prop: fiber comparison}. Together with the preceding determinant-and-projectivization comparison, the ordinary K\"unneth theorem therefore shows that the restriction of $\Psi_{\GL}$ to each fiber over $\M_{g,1}$ induces an isomorphism on rational cohomology through degree $R_{g,n}$. Applying \Cref{lem: range comparison} now gives
\[
\Psi_{\GL}^*\colon H^k(T_{g,n,d};\QQ)\longrightarrow H^k(M_B^{\GL}(\C_{g,1},n,d);\QQ)
\]
is an isomorphism for $k\leq R_{g,n}$.

We next compute the stable cohomology of $T_{g,n,d}$. Put $Y:=B\CC^*\times\BPGL_n$. Choose a topological abelian group model $B\CC^*\simeq K(\ZZ,2)$. Evaluation at the marked point then gives a $\Mod_{g,1}$-equivariant homotopy equivalence
\[
\Map(S_g,B\CC^*)_0
\simeq
\Map_*(S_g,B\CC^*)_0\times B\CC^*.
\]
Consequently,
\begin{equation}\label{eq: GL evaluation splitting}
\Sc_{g,1}(Y)_{(0,d)}\cong T_{g,n,d}\times B\CC^*.
\end{equation}
Let $\xi\in H^2\big(\Sc_{g,1}(Y)_{(0,d)};\QQ\big)$ be the pullback of the universal first Chern class from $B\CC^*$.

Combining \Cref{prop: marked point} with \eqref{eq: scanning algebra general}, and interpreting the formal variable $c_1$ in $K'_n$ as the universal first Chern class from the $B\CC^*$-factor of $Y$, gives a homomorphism
\[
\Sym^\bullet(K'_n)\otimes\QQ[\psi,\xi,c_2,\dots,c_n]
\longrightarrow
H^*\big(\Sc_{g,1}(Y)_{(0,d)};\QQ\big)
\]
which is an isomorphism in degrees at most $s(g)$ and sends the polynomial generator $\xi$ to the class defined above. Under \eqref{eq: GL evaluation splitting}, the K\"unneth formula identifies
\[
H^*\big(\Sc_{g,1}(Y)_{(0,d)};\QQ\big)
\cong H^*(T_{g,n,d};\QQ)\otimes\QQ[\xi].
\]
Quotienting the source and target by the ideal $(\xi)$ therefore gives a homomorphism
\begin{equation}\label{eq: GL target algebra}
\Xi\colon
\Sym^\bullet(K'_n)\otimes\QQ[\psi,c_2,\dots,c_n]
\longrightarrow H^*(T_{g,n,d};\QQ),
\end{equation}
which is again an isomorphism in degrees at most $s(g)$.

It remains to compare $\Xi$ with the tautological homomorphism of \eqref{eq: GLring}, which we denote by $\tau_{\GL}$. Let
\[
p_{\GL_n}\colon\mathcal C_{\GL_n}\longrightarrow M_B^{\GL}(\C_{g,1},n,d)
\]
be the universal curve, with universal marked section $\sigma_{\GL_n}$, and put $d_1:=c_1(\mathcal D_{\GL_n})$. By construction of $\Psi_{\GL}$ and naturality of Gysin maps,
\[
\Psi_{\GL}^*\Xi\big(\kappa_{a_0,a_1;\,a_2,\dots,a_n}\big)=(p_{\GL_n})_!\left(\psi_{\GL_n}^{a_0}d_1^{a_1}c_2(\mathcal P_{\GL_n})^{a_2}\cdots c_n(\mathcal P_{\GL_n})^{a_n}\right).
\]
On the marked generators, the same construction gives
\[
\Psi_{\GL}^*\Xi(\psi)=\sigma_{\GL_n}^*\psi_{\GL_n},
\qquad
\Psi_{\GL}^*\Xi(c_j)=\sigma_{\GL_n}^*c_j(\mathcal P_{\GL_n}),\qquad 2\leq j\leq n.
\]
These are precisely the images prescribed in \eqref{eq: GLring}, so $\tau_{\GL}=\Psi_{\GL}^*\circ\Xi$. The map $\Psi_{\GL}^*$ is an isomorphism through degree $R_{g,n}$, and $\Xi$ is an isomorphism through degree $s(g)$. Since $s(g)\leq R_{g,n}$, it follows that $\tau_{\GL}$ is an isomorphism in degrees at most $s(g)$.
\end{proof}

\appendix
\section*{Appendix: Stable intersection cohomology beyond the smooth case} 
\setcounter{section}{1}
\setcounter{theorem}{0}
\begin{center}
\textsc{Anne Larsen and Mirko Mauri}
\end{center}

The goal of this appendix is to extend Theorems \ref{thm:mainPGL} and \ref{thm:sl-stability-answer} beyond the smooth case $d \in (\ZZ/n\ZZ)^*$, to arbitrary degrees $d$. More precisely, we prove that for any $d \in \ZZ/n\ZZ$, the \textit{intersection cohomology} of the universal character variety stabilizes to a value independent of $d$. This is a relative version of the degree-independence of the stable intersection cohomology of Dolbeault moduli spaces on a single curve  \cite[Proposition~1.9]{MM}. 

\begin{theorem}\label{thm:mainappendix} 
    For any $d \in \ZZ/n\ZZ$ and $g\geq 2$, the intersection cohomology groups \[\IH^k(M^{\PGL}_B(\C_{g},n,d); \QQ), \IH^k(M^{\GL}_B(\C_{g,1},n,d); \QQ), \IH^k(M^{\SL}_B(\C_{g,1},n,d); \QQ)
    \] are independent of the degree $d$ for 
    $k < 4(g-1)(n-1)-2$. 
    In the overlap of these ranges, i.e., for $k \le R_{g,n}$ with $k < 4(g-1)(n-1)-2$, and when $d \in (\ZZ/n\ZZ)^*$, these intersection cohomology groups agree with the corresponding ordinary cohomology groups computed in \Cref{thm:mainPGL,thm:sl-stability-answer}.
\end{theorem}

\noindent The bound $4(g-1)(n-1)-2$ is optimal, except for $n=2$; see \S \ref{rmk:optimal}.

\vspace{0.1 cm}

\noindent\textbf{Notation.} Given $Z \subset X$ a closed and irreducible subvariety of a complex algebraic variety $X$, or more generally a closed embedding of integral Deligne--Mumford stacks, and $\mathcal{L}$ a local system on an open subset $U \subset Z$, we write $\IC_Z(\mathcal{L})$ for the intermediate extension of $\mathcal{L}$. We adopt the topologist's convention: $\IC_Z(\mathcal{L})$ is a bounded complex with trivial cohomology sheaves in negative degrees, and $\IC_Z(\mathcal{L})[\dim Z]$ is a perverse sheaf.
The intersection complex $\IC_X \coloneqq \IC_X( \QQ_X)$ is the intermediate extension of the constant sheaf, and its cohomology is the intersection cohomology with middle perversity and rational coefficients, denoted $\IH^*(X) \coloneqq H^*(X; \IC_X)$. The derived category of constructible complexes of sheaves on $X$ with rational coefficients is denoted $D^{b}_c(X, \QQ)
$. (Definitions can be found in \cite[\S 2]{dC-Migliorini}.)

\subsection{The direct image of the intersection complex along the Hitchin fibration}
\subsubsection{Spectral curves}
Fix $g\geq 2$. For any curve $[C] \in \M_g$, consider the affine spaces
$$\AA_0(C,n) \coloneqq \bigoplus_{i=2}^n H^0(C, \omega_C^{\otimes i}) \subset \bigoplus_{i=1}^n H^0(C, \omega_C^{\otimes i}) \eqqcolon \AA(C,n).$$
Each point $a = (a_i)_i \in \AA(C,n)$ defines the \textit{spectral curve} \[C^{\spec}_a \coloneqq V(t^n + (\pi^* a_1)t^{n-1}+ (\pi^* a_2) t^{n-2} + \ldots + \pi^* a_n) \subset T^*C\] of degree $n$ over $C$, where $t \in H^0(T^*C, \pi^* \omega_C)$ is the tautological section of the line bundle $\pi^* \omega_C$, i.e., the pullback of the canonical line bundle of $C$ along the projection $\pi \colon T^*C \to C$. This description allows us to define various loci in $\AA_0(C,n)$ and $\AA(C,n)$, in particular:

\begin{definition}
    The smooth (resp.\ elliptic) locus of $\AA_0(C,n)$ or $\AA(C,n)$ is the open subset over which the family of spectral curves is smooth (resp.\ integral).
\end{definition}

Recall the (proper) Hitchin morphisms
$$\chi^{\GL}\colon M^{\GL}_{\Dol}(C,n,d) \to \AA(C,n), \,\, \chi^{\SL}\colon M^{\SL}_{\Dol}(C,n,L) \to \AA_0(C,n),$$
$$\chi^{\PGL}\colon M^{\PGL}_{\Dol}(C,n,d) \to \AA_0(C,n),$$
which assign to each Higgs bundle $(E, \theta)$ the characteristic polynomial of its Higgs field $\theta$. %
The fibers of these morphisms are related to the spectral curves as follows:
\begin{itemize}
    \item $(\chi^{\GL})^{-1}(a)$ is isomorphic to the compactified Jacobian $\overline{\Jac}^d(C^{\spec}_a)$ parametrizing rank-one torsion-free sheaves on $C^{\spec}_a$ of degree $d+(n^2-n)(g-1)$, semistable with respect to the canonical polarization;
    \item $(\chi^{\SL})^{-1}(a)$ is isomorphic to $\overline{\Prym}^d(C^{\spec}_a/C)$, i.e., the fiber over $L$ of the norm map $
    \overline{\Jac}^d(C^{\spec}_a) 
    \to \Jac^d(C)$ given by $\mathcal I \mapsto \det(\pi_* \mathcal I)$;
    \item $(\chi^{\PGL})^{-1}(a)$ is isomorphic to the quotient $\overline{\Prym}^d(C^{\spec}_a/C)/\Jac^0(C)[n]$, where the action of $\Jac^0(C)[n]$ is by tensor product with the pullback bundle on $C_a^{\spec}$.
\end{itemize}
See \cite[\S 5.1]{Hitchin-stable} (over the smooth locus), \cite[Proposition 3.4]{BNR} (over the elliptic locus), and \cite{Schaub} and \cite[Proposition 6.1]{Hausel-Pauly} (in general).

The symmetries of the moduli spaces force tight relations between the direct images of the intersection complexes of $M^{G}_{\Dol}(C,n,d)$ along the maps $\chi^G$ for $G = \GL, \SL, \PGL$, which we explore in the following. We first compare $\chi^{\SL}$ and $\chi^{\PGL}$.

\subsubsection{$\SL$ vs $\PGL$}  %

The group $\Gamma \coloneqq \Jac^0(C)[n]$ acts on $M^{\SL}_{\Dol}(C,n,L)$ by tensor product, as follows: $L_{\gamma} \cdot (E, \theta)=(E \otimes L_\gamma, \theta)$ for any $L_{\gamma} \in \Gamma$. The Hitchin morphism $\chi^{\SL}$ is $\Gamma$-invariant. The direct image $R\chi_*^{\SL} \IC_{M^{\SL}_{\Dol}(C,n,L)}$ inherits a $\Gamma$-action and decomposes as the sum of the $\Gamma$-invariant and the $\Gamma$-variant (i.e., non-invariant)
parts: 
\begin{equation}\label{eq:SLsplit}
R\chi^{\SL}_* \IC_{M^{\SL}_{\Dol}(C,n,L)}= (R\chi^{\SL}_* \IC_{M^{\SL}_{\Dol}(C,n,L)})^\Gamma \oplus (R\chi^{\SL}_* \IC_{M^{\SL}_{\Dol}(C,n,L)})^{\mathrm{var}}.
\end{equation}
Since the $\Gamma$-quotient of $M^{\SL}_{\Dol}(C,n,L)$ is $M^{\PGL}_{\Dol}(C,n,d)$, Lemma \ref{lemma: IC finite quotient}
gives the following characterization of the invariant part. 
\begin{prop} \label{prop: SL fixed to PGL} There is an isomorphism in $D^{b}_c(\AA_0(C,n), \QQ)$
    $$(R\chi_*^{\SL} \IC_{M^{\SL}_{\Dol}(C,n,L)})^\Gamma \cong R\chi_*^{\PGL} \IC_{M^{\PGL}_{\Dol}(C,n,d)}.$$
\end{prop}

\begin{lemma} \label{lemma: IC finite quotient}
    Let $X$ be an irreducible complex variety acted on by a finite group $\Gamma$, and let $f\colon X \to X/\Gamma \eqqcolon Y$ be the geometric quotient. Then \[(f_* \IC_X)^\Gamma \cong \IC_Y.\]
\end{lemma}

\begin{proof}
    Let $U \subset Y$ be a dense smooth open subset over which $f$ is \'etale. Since $f$ is finite, $f_* \IC_X$ is an intersection complex, as it satisfies the support and cosupport conditions for intersection complexes \cite[(12), (13)]{dC-Migliorini}. This means that \begin{equation}\label{eq:intcoh}
        f_* \IC_X \cong \IC_Y(f_* \QQ_{f^{-1}(U)}).
    \end{equation} 
    The action of $\Gamma$ on $X$ induces $\Gamma$-equivariant structures on both sides of \eqref{eq:intcoh}, and the intermediate-extension isomorphism is $\Gamma$-equivariant.
We conclude that
\[
(f_* \IC_X)^\Gamma = \IC_Y\big((f_* \QQ_{f^{-1}(U)})^\Gamma\big) = \IC_Y(\QQ_{U}). \qedhere
\]
\end{proof}
The invariant part has no proper direct summand on the elliptic locus.
\begin{prop}\label{prop:fullsupportPGLell}
The complex $R\chi_*^{\PGL} \IC_{M^{\PGL}_{\Dol}(C,n,d)}$ has full support on the elliptic locus.
In particular, the codimension of each proper support is $\geq 2(g-1)(n-1)-1$. 
\end{prop}
\begin{proof}
By \cite[\S 7]{Ngo}, a Ng\^{o} fibration with integral fibers (in particular, the Hitchin fibration over the elliptic locus) has full support. Alternatively, see \cite[Theorem 3.13]{MM} together with Proposition \ref{prop: GL and SL}, or \cite{dCHM2021} 
together with Proposition \ref{prop: GL and SL} and Proposition \ref{prop:chiindependence}.
\end{proof}

\begin{prop}\label{prop:chiindependence}
Let $e \in (\ZZ/n\ZZ)^*$ and $L' \in \Pic^{e}(C)$. Then \[R\chi_*^{\GL} \IC_{M^{\GL}_{\Dol}(C,n,d)}, \quad R\chi_*^{\PGL} \IC_{M^{\PGL}_{\Dol}(C,n,d)}, \quad R\chi^{\SL}_* \IC_{M^{\SL}_{\Dol}(C,n,L)}\] are respectively direct summands of \[R\chi_*^{\GL} \IC_{M^{\GL}_{\Dol}(C,n,e)}, \quad R\chi_*^{\PGL} \IC_{M^{\PGL}_{\Dol}(C,n,e)}, \quad R\chi^{\SL}_* \IC_{M^{\SL}_{\Dol}(C,n,L')}.\]
\end{prop}
\begin{proof}
This follows from the degree-independence of BPS sheaves \cite{Kinjo-Koseki} %
in the case of $G=\GL$ by \cite[Proposition 3.9]{MM}, thus also for $G=\PGL$ by Proposition \ref{prop: GL and SL}. The same argument works for $G=\SL$ mutatis mutandis; cf. \cite[Lemma 8.4]{Padurariu-Toda}. %
Note that $R\chi^{\SL}_* \IC_{M^{\SL}_{\Dol}(C,n,L)}$ is independent of $L \in \Pic^d(C)$.
\end{proof}

The variant part instead is supported on the endoscopic locus, which we now define.
\begin{definition} 
    The endoscopic locus $\AA_0(C,n)^{\nne}$ is the locus of points $a\in \AA_0(C,n)$ such that the group scheme $\Prym^{0}(C^{\spec}_a/C)$ is disconnected, where $\Prym^{0}(C^{\spec}_a/C)$ is the kernel of the norm map $\Jac^{0}(C^{\spec}_a) \to \Jac^0(C)$ given by $\mathcal{L} \mapsto \det(\pi_*\mathcal{L}) \otimes \det(\pi_*\mathcal{O}_{C^{\spec}_a})^{\vee}$. 
\end{definition}

\begin{prop} \label{prop: SL var}
     The complex $(R\chi^{\SL}_* \IC_{M^{\SL}_{\Dol}(C,n,L)})^{\mathrm{var}}$ is supported on the endoscopic locus. In particular, its support has codimension at least $\rho_{g,n}/2$ in $\AA_0(C,n)$, with $\rho_{g,n}$ as in \eqref{eq:rho-range}.
\end{prop} 
\begin{proof} By Proposition \ref{prop:chiindependence}, it suffices to consider the case $\deg(L) \in (\ZZ/n\ZZ)^*$, so that $M^{\SL}_{\Dol}(C,n,L)$ is smooth, and $\IC_{M^{\SL}_{\Dol}(C,n,L)}=\QQ_{M^{\SL}_{\Dol}(C,n,L)}$. Now, let $\IC_Z(\mathcal{L})[-k]$ be a simple proper direct summand of $(R\chi^{\SL}_* \QQ_{M^{\SL}_{\Dol}(C,n,L)})^{\mathrm{var}}$, so that in particular $\Gamma$ acts on $\mc L$ by some nontrivial character. Suppose by contradiction that $Z$ does not lie in the endoscopic locus, and choose a general point $a \in Z$. Then $\mathcal{L}_a$ is a $\Gamma$-equivariant direct summand of $(R^{k} \chi^{\SL}_{*} \QQ_{M^{\SL}_{\Dol}(C,n,L)})_a =H^k((\chi^{\SL})^{-1}(a); \QQ)$.
Since $a$ does not lie in the endoscopic locus, the action of $\Gamma$ on the fiber $(\chi^{\SL})^{-1}(a)$ extends to the action of the connected group $\Prym^{0}(C^{\spec}_a/C)$, so
$\Gamma$ acts trivially on the fiberwise cohomology groups $H^*((\chi^{\SL})^{-1}(a); \QQ)$. However, since $\Gamma$ acts on $\mc L_a$ by a nontrivial character, we see that $\mathcal{L}=0$, i.e., the variant part of $R\chi^{\SL}_* \IC_{M^{\SL}_{\Dol}(C,n,L)}$ is supported on the endoscopic locus. %
Finally, the codimension of the endoscopic locus is $\rho_{g,n}/2$ by \cite[Corollary 1.3, Lemma 7.1]{Hausel-Pauly} and \eqref{eq:rho-range}. %
\end{proof}

\begin{remark}
The reduction to the smooth case in the proof of Proposition~\ref{prop: SL var} %
requires the degree-independence of BPS sheaves (Proposition~\ref{prop:chiindependence}). Alternatively, one may observe that the integer $k$ in the previous proof can be taken to be $2c$, where $c$ is the relative dimension of $\chi^{\SL}$, by Ng\^{o} freeness and $\delta$-regularity; see \cite[Proposition 1.5]{chi-indep} and \cite[Theorem A]{dCFHM}. In top degree, the direct images of the intersection complex and constant sheaf coincide by \cite[Proposition~3.3]{MM}, i.e., $R^{2c} \chi^{\SL}_{*} \IC_{M^{\SL}_{\Dol}(C,n,L)} \cong R^{2c} \chi^{\SL}_{*} \QQ_{M^{\SL}_{\Dol}(C,n,L)}$, which makes the proof work in any degree.  For our purposes, deep arguments using the degree-independence or Ng\^{o} freeness can be avoided: it suffices to prove the proposition over the elliptic locus, over which $M^{\SL}_{\Dol}(C,n,L)$ is smooth in any degree.
The codimension estimate of the complement of the nonendoscopic elliptic locus is worse, namely $\geq 2(g-1)(n-1)-1$, but sufficient for  \Cref{thm:mainappendix}.
\end{remark}

 \subsubsection{$\GL$ vs $\PGL$}
Recall that $H^0(C, \omega_C)$ acts on $\AA(C,n)$ via the translation action of $H^0(C, \omega_C) \xrightarrow{\pi^*} H^0(T^*C, \pi^*\omega_C)$ on the tautological section $t \in H^0(T^*C, \pi^* \omega_C)$.
Consider the quotient map $q \colon \AA(C,n) \cong H^0(C, \omega_C) \oplus \AA_0(C,n)\to \AA_0(C,n)$ given explicitly by
$$p(t) = t^n + \sum^{n}_{i=1} \pi^*a_{i}t^{n-i} \mapsto p(t-\pi^*a_1/n)\eqqcolon t^n + \sum^{n}_{i=2} \pi^*b_{i}t^{n-i}.$$ %

\begin{prop}\label{prop: GL and SL}
    There is an isomorphism in $D^{b}_c(\AA(C,n), \QQ)$
   \[
    R\chi_*^{\GL} \IC_{M^{\GL}_{\Dol}(C,n,d)} \cong H^*(\Pic^0(C); \QQ)\otimes q^* R\chi^{\PGL}_* \IC_{M^{\PGL}_{\Dol}(C,n,d)}.
 \]
 
\end{prop}
\begin{proof}
Recall that there is a commutative diagram 
\[
\xymatrix{
M^{\GL}_{\Dol}(C,1,0) \times M^{\SL}_{\Dol}(C,n,L) \ar[r]^-{/\, \Gamma} \ar[d]_{(\chi^{\mathbb{G}_m}, \chi^{\SL})} & M^{\GL}_{\Dol}(C,n,d) \ar[d]^{\chi^{\GL}}\\
H^0(C, \omega_C) \times \AA_0(C,n) \ar[r]^-{\cong} & \AA(C,n)
}
\]
where the top horizontal arrow is given by
$$((E, \theta), (F, \psi)) \mapsto \left(E \otimes F, \id_E \otimes \psi + \theta \otimes \id_F\right)$$
and the bottom horizontal arrow by $(s, p(t)) \mapsto p(t+\pi^*s)$. In fact, the top horizontal arrow is the quotient by the free $\Gamma$-action given by
$$L_{\gamma} \cdot ((E, \theta), (F, \psi)) = ((E \otimes L_{\gamma}, \theta), (F \otimes L_{\gamma}^{-1}, \psi))$$ for any $L_{\gamma}\in \Gamma$.
By Lemma \ref{lemma: IC finite quotient} %
and the K\"{u}nneth formula in intersection cohomology \cite[Proposition 2]{IH-Kunneth}, we obtain
\begin{align*}
R \chi^{\GL}_* \IC_{M^{\GL}_{\Dol}(C,n,d)} 
& \cong (R (\chi^{\mathbb{G}_m}, \chi^{\SL})_*\IC_{M^{\GL}_{\Dol}(C,1,0) \times M^{\SL}_{\Dol}(C,n,L)})^{\Gamma}\\
& \cong (
R \chi^{\mathbb{G}_m}_* 
\IC_{M^{\GL}_{\Dol}(C,1,0)}  
\boxtimes
R \chi^{\SL}_*\IC_{M^{\SL}_{\Dol}(C,n,L)}
)^{\Gamma}.
\end{align*}
Since the Hitchin map $\chi^{\mathbb{G}_m} \colon M^{\GL}_{\Dol}(C,1,0) \to H^{0}(C, \omega_C)$ is isomorphic to the projection $\Pic^0(C) \times H^{0}(C, \omega_C) \to H^{0}(C, \omega_C)$, we have
\[R \chi^{\mathbb{G}_m}_* 
\IC_{M^{\GL}_{\Dol}(C,1,0)} \cong \QQ_{H^0(C, \omega_C)} \otimes H^*(\Pic^0(C); \QQ).\] 
Further, since the translation action of $\Gamma$ on $\Pic^0(C)$ induces the trivial action in cohomology, we obtain 
\begin{align*}
(
R \chi^{\mathbb{G}_m}_* 
& \IC_{M^{\GL}_{\Dol}(C,1,0)}  
 \boxtimes
R \chi^{\SL}_*\IC_{M^{\SL}_{\Dol}(C,n,L)}
)^{\Gamma}\\
& \cong \big(\QQ_{H^0(C, \omega_C)} \otimes H^*(\Pic^0(C); \QQ) \boxtimes  R \chi^{\SL}_*\IC_{M^{\SL}_{\Dol}(C,n,L)}\big)^{\Gamma}\\
& \cong H^*(\Pic^0(C); \QQ) \otimes  q^* (R \chi^{\SL}_*\IC_{M^{\SL}_{\Dol}(C,n,L)})^{\Gamma}\\
& \cong H^*(\Pic^0(C); \QQ) \otimes q^* R \chi^{\PGL}_*\IC_{M^{\PGL}_{\Dol}(C,n,d)},
\end{align*}
where the last isomorphism follows from Proposition \ref{prop: SL fixed to PGL}.
\end{proof}
\subsubsection{Shifts of summands and codimension of supports}
\begin{lemma}\label{propboundcoeff}
Let $f \colon X \to Y$ be a proper
morphism of algebraic varieties satisfying $R^{> 2c}f_* \IC_X =0$ for $c := \dim X - \dim Y$. For all direct summands $\IC_{Z}(\mathcal{L})[-k]$ of $Rf_* \IC_X$ with $Z \subsetneq Y$, we have 
\[
k\geq 2 \codim Z.
\]
\end{lemma}
\begin{proof}
If $\IC_{Z}(\mathcal{L})[-k]$ is a direct summand of $Rf_* \IC_X$, then $\IC_{Z}(\mathcal{L})[\dim Z][c-k+\codim Z]$ is a direct summand of $Rf_* \IC_X[\dim X]$. By the properness of $f$ and since $\IC_X[\dim X]$ is self-dual, Poincar\'{e}-Verdier duality gives 
\[{}^{\mathfrak{p}} \mathcal{H}^{c-k+\codim Z}(Rf_* \IC_X[\dim X]) \cong {}^{\mathfrak{p}} \mathcal{H}^{-(c-k+\codim Z)}(Rf_* \IC_X[\dim X])^\vee,\]
which implies that $\IC_{Z}(\mathcal{L}^\vee)[-k-2(c-k+\codim Z)]$ is also a direct summand of $Rf_* \IC_X$. Since $R^{> 2c}f_* \IC_X =0$, we must have $k+2(c-k+\codim Z) \leq 2c$, i.e., $k\geq 2 \codim Z$. 
\end{proof}

\begin{remark}\label{rmkcohomologicalbound}
The cohomological bound $R^{> 2c}f_* \IC_X =0$ holds if $X$ admits a Whitney stratification $X = \bigsqcup X_{i}$ such that for any $y \in Y$
\[
\codim_{X}(X_{i}) = 2 \codim_{f^{-1}(y)}(X_{i}\cap f^{-1}(y));
\]
 see \cite[Proposition 3.3]{MM}. This is the case if $f$ is a Lagrangian fibration from a possibly singular symplectic variety \cite[Proposition 3.6]{MM}. In particular, it holds for the Hitchin fibrations $\chi^G$ above, for $G=\GL,\SL,\PGL$, by \cite[Proposition 1.5]{dCFHM}. 
\end{remark}

\subsection{Degree-independence of the stable intersection cohomology of the universal character variety}
\subsubsection{The case of $\PGL$}
Recall the relative Hitchin fibration \[\chi^{\PGL}\colon M^{\PGL}_{\Dol}(\C_g,n,d) \to \AA_0(\C_g,n),\] where $\AA_0(\C_g,n) \to \M_g$ is the affine bundle with fiber $\AA_0(C,n)$ over $[C]$. 
Consider the universal spectral curve $\mathscr C^{\spec} \to \AA_0(\C_g,n)$ over $\C_g \times_{\M_g} \AA_0(\C_g,n)$. Let $\AA_0(\C_g,n)^\circ$ be the locus over which $\mathscr C^{\spec}$ is smooth.

\begin{prop} \label{prop: torsor}
   The complex $R\chi^{\PGL}_* \IC_{M^{\PGL}_{\Dol}(\C_g,n,d)}|_{\AA_0(\C_g, n)^\circ}$ has full support and is independent of $d$.
\end{prop}

\begin{proof}
   Let $\Pic^0(\C_g \times_{\M_g} \AA_0(\C_g,n)/\AA_0(\C_g,n))$ (resp.\ $\Pic^0(\mathscr C^{\spec}/\AA_0(\C_g,n))$) be the Picard functor of the pullback of the universal curve $\mathcal{C}_g \to \M_g$ (resp.\ of the universal spectral curve) over $\AA_0(\C_g,n)$. The kernel of the norm map 
   \begin{align*}
       \Pic^0(\mathscr C^{\spec}/\AA_0(\C_g,n)) & \to \Pic^0(\C_g \times_{\M_g} \AA_0(\C_g,n)/\AA_0(\C_g,n))
    \\
    \mathcal{L} & \mapsto \det(\pi_*\mathcal{L}) \otimes \det(\pi_*\mathcal{O}_{\mathscr C^{\spec}})^{\vee}
   \end{align*}
is the group stack denoted $\Prym^0(\mathscr C^{\spec}/\C_g \times_{\M_g} \AA_0(\C_g,n))$.
   Define the group stack
    \begin{equation}\label{groupscheme}
        \mathscr P \coloneqq \Prym^0(\mathscr C^{\spec}/\C_g \times_{\M_g} \AA_0(\C_g,n))/\Gamma(\C_g),
    \end{equation}
    where $\Gamma(\C_g)$ is the $n$-th torsion of the universal Picard functor over $\M_g$, pulled back over $\AA_0(\C_g,n)$, i.e., $\Pic^0(\C_g/\M_g)[n] \times_{\M_g} \AA_0(\C_g,n)$. Note that the restriction of $M^{\PGL}_{\Dol}(\C_g,n,d)$ to $\AA_0(\C_g,n)^\circ$ is a (smooth) $\mathscr P$-torsor over $\AA_0(\C_g,n)^\circ$. Then over $\AA_0(\C_g,n)^\circ$ the intersection complex $\IC_{M^{\PGL}_{\Dol}(\C_g,n,d)}$ is the constant sheaf $\QQ$, and by \cite[Lemma 1.3.5]{PW} we see that the restrictions of $R\chi^{\PGL}_* \IC_{M^{\PGL}_{\Dol}(\C_g,n,d)}$ and $R\chi_* \QQ_{\mathscr P}$ to $\AA_0(\C_g,n)^\circ$ are isomorphic. The latter is independent of $d$.   
\end{proof}
    Therefore, the only part of $R\chi^{\PGL}_* \IC_{M^{\PGL}_{\Dol}(\C_g,n,d)}$ that depends on $d$ is supported on the complement of $\AA_0(\C_g,n)^\circ$. Direct summands in the next proposition are taken in $D^b_c(\AA_0(\C_g,n), \QQ)$.

\begin{prop} \label{prop: PGL d dependent}
If $\IC_{Z}(\mathcal{L})[-k]$ is a direct summand of $R\chi^{\PGL}_* \IC_{M^{\PGL}_{\Dol}(\C_g,n,d)}$ for some $Z \subsetneq \AA_0(\C_g,n)$, then $k \geq 4(g-1)(n-1)-2$.
\end{prop}

\begin{proof}
The analogous statement holds fiberwise for any $\AA_0(C,n)$, but to prove it over $\AA_0(\C_g,n)$, we do the following: Choose $[C] \in \M_g$ such that $Z \cap \AA_0(C, n)\neq \emptyset$. By construction, the relative Betti moduli space $M^{\PGL}_{B}(\C_g,n,d) \to \M_g$ is a locally trivial family, and so the relative Dolbeault moduli space $M^{\PGL}_{\Dol}(\C_g,n,d) \to \M_g$ is too by Simpson's relative nonabelian Hodge correspondence \cite[Theorem 9.11, Lemma 9.14]{Simpson-reps2}. The local triviality implies that
\[\IC_{M^{\PGL}_{\Dol}(\C_g,n,d)}|_{M^{\PGL}_{\Dol}(C,n,d)} \cong \IC_{M^{\PGL}_{\Dol}(C,n,d)}.\]
Therefore, $\IC_{Z}(\mathcal{L})[-k]|_{\AA_0(C,n)}$ is a direct summand of 
\[
R\chi^{\PGL}_* \IC_{M^{\PGL}_{\Dol}(\C_g,n,d)}|_{\AA_{0}(C, n)} \cong R\chi^{\PGL}_* \IC_{M^{\PGL}_{\Dol}(C,n,d)} \in D^b_c(\AA_0(C,n)).
\]
Note that $Z \cap \AA_0(C, n) \subsetneq \AA_0(C, n)$, since $Z$ avoids $\AA_0(\C_g, n)^\circ$ by Proposition \ref{prop: torsor}, and $\AA_0(\C_g, n)^\circ$ is dense in each $\AA_0(C, n)$. By Proposition \ref{prop:fullsupportPGLell}, \[\codim_{\AA_0(C, n)}(Z \cap \AA_0(C, n))\geq 2(g-1)(n-1)-1\] for all $[C] \in \M_g$ such that $Z \cap \AA_0(C, n)\neq \emptyset$. In particular, we have
\[\codim_{\AA_0(\mathcal{C}_g, n)}(Z)\geq 2(g-1)(n-1)-1\] and so the result follows from Lemma \ref{propboundcoeff} and Remark \ref{rmkcohomologicalbound}.
\end{proof}

\begin{prop} \label{prop: PGL result} The intersection cohomology
    $$\IH^{k}(M^{\PGL}_{\Dol}(\C_g,n,d))
    \cong \IH^{k}(M^{\PGL}_{B}(\C_g,n,d))$$ is independent of $d$ for $k< 4(g-1)(n-1)-2$.
\end{prop}

\begin{proof}
    The direct image $R\chi_*^{\PGL} \IC_{M_{\Dol}^{\PGL}(\C_g,n,d)}$ decomposes into a sum of simple shifted perverse sheaves fully supported on $\AA_0(\C_g,n)$ and a sum of those supported on proper subvarieties of $\AA_0(\C_g,n)$. By Proposition \ref{prop: torsor} the first term is independent of $d$. By Proposition \ref{prop: PGL d dependent} the cohomology of the second term vanishes in degree $<  4(g-1)(n-1)-2$. So $\IH^k(M^{\PGL}_{\Dol}(\C_g,n,d); \QQ)$ is independent of $d$ for $k <  4(g-1)(n-1)-2$.
\end{proof}

\subsubsection{From $\PGL$ to $\GL$ and $\SL$}
Lacking a relative  nonabelian Hodge correspondence for degree $d \ne 0$ and $G = \GL$ or $\SL$, we reduce \Cref{thm:mainappendix} for these groups to the corresponding statement for $\PGL$ via the fiberwise nonabelian Hodge correspondence. Let $\pi_G\colon M^G_B(\C_{g,1},n,d) \to \M_{g,1}$ be the projection map, and set $\mathcal{K}^{G}_{B} \coloneqq R\pi_{G, *} \IC_{M^G_{B}(\C_{g,1},n,d)}$.

\begin{prop} \label{prop: SL result}
With $\rho_{g,n}$ as in \eqref{eq:rho-range}, set $N:=4(g-1)(n-1)-2$. There is an isomorphism
$$\IH^k(M^{\SL}_B(\C_{g,1},n,d); \QQ) \cong \IH^k(M^{\PGL}_B(\C_{g,1},n,d); \QQ) $$
    for $k < \rho_{g,n}$. Furthermore, these groups are independent of $d$ for $k < N$. %
\end{prop}

\begin{proof}
Since the morphism $f\colon M^{\SL}_B(\C_{g,1},n,d) \to M_B^{\PGL}(\C_{g,1},n,d)$ is a quotient by the relative group object $\widetilde\Gamma$ of Proposition~\ref{prop: relative B relations}, we write, as in \eqref{eq:SLsplit} and Proposition \ref{prop: SL fixed to PGL}, that
\begin{equation}\label{eq:splitfinalSL}
    \mathcal{K}^{\SL}_{B} \cong (\mathcal{K}^{\SL}_{B})^{\widetilde\Gamma} \oplus (\mathcal{K}^{\SL}_{B})^{\mathrm{var}} \cong \mathcal{K}^{\PGL}_{B} \oplus (\mathcal{K}^{\SL}_{B})^{\mathrm{var}}.
\end{equation}
The local triviality of  $\pi_G$ implies that \[(\mathcal{K}^{\SL}_{B})^{\mathrm{var}}|_{[C]} \cong \IH^*(M^{\SL}_{B}(S_{g,1}, n, d))^{\mathrm{var}}.\]
By the fiberwise ($\Gamma$-equivariant) nonabelian Hodge correspondence for $\SL$ \cite[Theorem 4.3]{Hausel-handbook} (cf.~\cite[Remark 5 after Theorem 1]{H2024}), there is an isomorphism
\[
\IH^k(M^{\SL}_{B}(S_{g,1}, n, d))^{\mathrm{var}} \cong \IH^k(M^{\SL}_{\Dol}(C, n, L))^{\mathrm{var}}.
\]
This group vanishes for $k< \rho_{g,n}$ by Proposition \ref{prop: SL var} and Lemma \ref{propboundcoeff}. Hence, taking global cohomology of \eqref{eq:splitfinalSL}, we conclude that $\IH^k(M^{\SL}_B(\C_{g,1},n,d); \QQ) \cong \IH^k(M^{\PGL}_B(\C_{g,1},n,d); \QQ)$
    for $k < \rho_{g,n}$. As in Proposition \ref{prop: PGL result}, these intersection cohomology groups are independent of $d$ for $k < N$.
\end{proof}

\begin{prop}\label{main resultGL} The intersection cohomology  $\IH^k(M^{\GL}_B(\C_{g,1},n,d); \QQ)$ is independent of $d$ for
    $k < 4(g-1)(n-1)-2$.
\end{prop}

\begin{proof}
By the argument of \eqref{eq:RpiGL} and Lemma \ref{lemma: IC finite quotient}, we have 
$\mathcal{K}^{\GL}_{B} \cong \mathcal{K}^{\PGL}_{B} \otimes \mathcal{K}^{\mathbb{G}_m}_{B}$.
Let $V_{g}$ be the shifted local system on $\mathcal{M}_{g,1}$ corresponding to the standard symplectic representation of the mapping class group in cohomological degree 1 in $D^b_c(\M_{g,1}, \QQ)$. Since  $\mathcal{K}^{\mathbb{G}_m}_{B} \cong \bigoplus_i \Lambda^i V_g$ is independent of $d$, the result follows from Proposition \ref{prop: PGL result}.
\end{proof}

\begin{proof}[Proof of Theorem \ref{thm:mainappendix}]
    This is the content of Propositions \ref{prop: PGL result} (for $\PGL$), \ref{prop: SL result} (for $\SL$), and \ref{main resultGL} (for $\GL$).
\end{proof}

\subsubsection{Sharpness of the bounds}\label{rmk:optimal} We discuss the sharpness of Propositions \ref{prop: PGL d dependent} and \ref{prop: PGL result}.
\begin{prop}\label{prop:boundI}
The bound $N\coloneqq 4(g-1)(n-1)-2$ is sharp for Proposition \ref{prop: PGL d dependent} for arbitrary values of the pair $(n,g)$.
\end{prop}
\begin{proof} Consider the inclusion $i \colon Z^{\circ} \hookrightarrow \AA_0(\mathcal{C}_g,n)$ of the locus of nodal spectral curves which are unions of two irreducible components, one of which is isomorphic to the underlying curve $C$, i.e., such that the two components of the spectral curve have arithmetic genus $g$ and $(n-1)^2(g-1)+1$ respectively, and let $Z$ be its closure. By \cite[Theorem 1.1]{MM}, there is a direct summand of $R \chi^{\PGL}_* \QQ_{M^{\PGL}_{\Dol}(\C_g,n,1)}$, which is not a summand of $R \chi^{\PGL}_* \IC_{M^{\PGL}_{\Dol}(\C_g,n,0)}$, of the form $i_* \mathcal{L}[-N]$ where $\mathcal{L}$ is a rank one local system on $Z^{\circ}$. 
\end{proof}
\begin{prop}\label{prop:boundII}
The bound $N\coloneqq 4(g-1)(n-1)-2$ in Proposition \ref{prop: PGL result} is sharp, except for $n=2$. 
\end{prop}
\begin{proof} 
 Let $i_* \mathcal{L}[-N]$ be the direct summand of $R \chi^{\PGL}_* \QQ_{M^{\PGL}_{\Dol}(\C_g,n,1)}$ in the proof of Proposition \ref{prop:boundI}. By relative hard Lefschetz, there is a corresponding summand $\mathcal{L}$ of the sheaf $R^{\textrm{top}}\chi_*^{\PGL} \QQ_{M^{\PGL}_{\Dol}(\C_g,n,1)}|_{Z^\circ}$. The fiber of this sheaf at the point $a \in Z^\circ$ is generated by the fundamental classes of the irreducible components of $(\chi^{\PGL})^{-1}(a)$. 
The fiber $(\chi^{\PGL})^{-1}(a)$ is a compactified Jacobian of $C^{\spec}_{a}$ whose irreducible components are indexed by the degrees of the restriction of a stable line bundle on $C^{\spec}_a$ to its irreducible components; cf.~\cite[\S 4]{dCHM2021} and \cite[Proposition 7.6.(ii)]{MMP2023}. As described in the proof of Proposition \ref{prop:boundI}, the spectral curve $C^{\spec}_a$ has two irreducible components, which have (crucially) different genus when $n \neq 2$. 
In particular, they cannot be exchanged by the monodromy action, and thus we conclude that $R^{\mathrm{top}}\chi_*^{\PGL} \QQ_{M^{\PGL}_{\Dol}(\C_g,n,1)}|_{Z^\circ}$ is a trivial local system. The same must then be true of its summand $\mathcal{L}$.

Now, as $\mathcal{L}$ is of rank one, for $n \neq 2$, the contribution of $i_* \mathcal{L}[-N] = i_* \QQ_{Z^\circ}[-N]$ to $H^N(M^{\PGL}_{\Dol}(\C_g,n,1))$ is $H^0(Z; i_* \QQ_{Z^\circ}) = \QQ$. Since this summand does not occur in the decomposition for $d=0$, we conclude that for $n \neq 2$
$$\dim H^{N}(M^{\PGL}_{\Dol}(\C_g,n,1)) \geq \dim \IH^{N}(M^{\PGL}_{\Dol}(\C_g,n,0))+1;$$
in particular the two groups differ, so the bound $N$ is sharp for $n \neq 2$.

On the other hand, if $n=2$, the monodromy of $\mathcal{L}|_{\AA_{0}(C, n) \cap Z^{\circ}}$ is no longer trivial. The locus $\AA_0(C,2) \cap Z^\circ$ is the image of the squaring map $H^{0}(C, \omega_C)\setminus 0 \to H^{0}(C, \omega^{\otimes 2}_C) \setminus 0$; since this is a covering map from the simply connected space $H^0(C, \omega_C) \setminus 0$, we conclude that $\pi_1(\AA_0(C,2) \cap Z^\circ)=\ZZ/2\ZZ$. The monodromy along this loop acts on the fiber of $\mathcal{L}$ by $-\id$; see \cite[Corollary 6.20]{dCFHM} or \cite[Theorem 1.6]{MMP2023}. Hence, if $n=2$, we have $H^0(Z; i_* \mathcal{L}) = 0$ and
$\dim H^{N}(M^{\PGL}_{\Dol}(\C_g,2,1)) = \dim \IH^{N}(M^{\PGL}_{\Dol}(\C_g,2,0)).$
\end{proof} 

\noindent\textbf{Acknowledgment.} Mirko Mauri was supported by Université Paris Cité and Sorbonne Université, CNRS, IMJ-PRG, F-75013 Paris, France.

\bibliographystyle{amsalpha}
\bibliography{references}
\end{document}